\documentclass{amsproc}
\usepackage{amssymb}
\usepackage{amsfonts}

\theoremstyle{plain}

\newtheorem{definition}{Definition}

\newtheorem{lemma}{Lemma}

\newtheorem{proposition}{Proposition}

\newtheorem{theorem}{Theorem}
\numberwithin{equation}{section}
\input{tcilatex}

\begin{document}
\title[Short Title]{Functionally recursive rings of matrices-Two examples}
\author{Said N. Sidki}
\address{Departamento de Matem\'{a}tica, Universidade de Bras\'{\i}lia,
Brasilia, DF 70910-900, Brazil}
\email{sidki@mat.unb.br}
\thanks{The author thanks Roland Bacher, Laurent Bartholdi and Pierre de la
Harpe for stimulating conversations on this topic. The author acknowledges
support from the Brazilian Conselho Nacional de Pesquisa and from FAPDF. }
\date{Dec. 20, 2008}

\begin{abstract}
We define the notions of finite-state and functionally recursive matrices
and their growth. We also introduce two rings generated by functionally
recursive matrices. The first is isomorphic to the $2$-generated free ring.
The second is a $2$-generated monomial ring such that the multiplicative
semigroup of monomials in the generators is nil of degree $5$ and the ring
has Gelfand Kirillov dimension $1+\frac{\log \left( 2\right) }{\log \left(
\alpha \right) }$ where $\alpha =\frac{1+\sqrt{5}}{2}$.
\end{abstract}

\maketitle

\section{Introduction}

Infinite dimensional matrices which are recursively defined and rings
generated by them have received attention stimulated to some degree by works
on representations of the group algebras of the infinite torsion groups of
Grigorchuk and of Gupta-Sidki \cite{sid97}, \cite{vieira01}, \cite%
{bar-grig02}, \cite{nek04}. The most recent formulations of notions of
matrix recursion are due to L. Bartholdi in \cite{bar06} and in the software 
\cite{bar08} with applications to group algebras and by R. Bacher in \cite%
{bach08} with applications to Dirichlet characters.

In the first part of this paper, we consider vector spaces having for base a
finitely generated free monoid and we consider linear transformations which
leave invariant the span of monoid elements of equal length. Following the
model of tree automorphisms, we define or re-define the notions of
finite-state and functionally recursive matrix transformations. We also
extend the notion of growth of tree automorphisms to matrix transformations.

In the second and longest part of the paper we introduce two recursive rings
as case studies. With the proper definition of recursive matrices, we prove
the following results.

\begin{theorem}
Let $R_{1}$ be the ring generated by $s=\left( 
\begin{array}{cc}
1 & 0 \\ 
0 & 2s%
\end{array}%
\right) ,t=\left( 
\begin{array}{cc}
0 & 2s \\ 
0 & 2t%
\end{array}%
\right) $ over $\mathbb{Z}$. Then $R_{1}$ is isomorphic to the free $2$%
-generated ring.
\end{theorem}

The matrices $s,t$ in $R_{1}$ are functionally recursive of exponential
growth. Even though they are not finite-state, the linear spans of their
states have finite rank.

\begin{theorem}
Let $R_{2}$ be the ring generated by $s=\left( 
\begin{array}{cc}
0 & 0 \\ 
1 & 0%
\end{array}%
\right) ,t=\left( 
\begin{array}{cc}
0 & t \\ 
0 & s%
\end{array}%
\right) $ over $\mathbb{Z}$. \newline
(i) Let $\mu \left( R_{2}\right) $ be the multiplicative semigroup generated
by $s,t$. Then $\mu \left( R_{2}\right) $ is nil of degree $5$. \newline
(ii) Let $\lambda $ be the endomorphism of the free ring $F=\left\langle
\sigma ,\tau \right\rangle $ defined by $\lambda :\sigma \rightarrow \tau
,\tau \rightarrow \sigma \tau $. Then $R_{2}$ affords the monomial
presentation $\left\langle \sigma ,\tau \mid u_{i}=0\text{ }\left( i\geq
1\right) \right\rangle $ where 
\begin{eqnarray*}
u_{1} &=&\sigma ^{2}, \\
u_{2i} &=&\tau .\left( u_{2i-1}\right) ^{\lambda }, \\
u_{2i+1} &=&\left( u_{2i}\right) ^{\lambda }\text{ }\left( i\geq 1\right) 
\text{.}
\end{eqnarray*}%
\newline
(iii) The Gelfand-Kirillov dimension of $R_{2}$ is $1+\frac{\log \left(
2\right) }{\log \left( \alpha \right) }$ where $\alpha =\frac{1+\sqrt{5}}{2}$%
.
\end{theorem}

The matrices $s,t$ in $R_{2}$ are finite-state and have bounded growth. More
detailed information about the structure of $R_{2}$ and its quotients can be
found in the text.

\section{Tree indexed linear spaces and linear transformations}

Let $Y=\left\{ y_{1,}y_{2,}...,y_{m}\right\} $ be a finite alphabet with $%
m\geq 2$ elements, $S$ the set of finite sequences on the alphabet $Y$ (or
words from the free monoid freely generated by $Y$), $S_{n}$ be the subset
of $S$ of sequences of length $n\geq 0$ and denote the empty sequence by $%
\phi $. Then, 
\begin{equation*}
S=\left\{ \phi \right\} \cup \left\{ yS\mid y\in Y\right\} =\cup _{n\geq
0}S_{n}\text{.}
\end{equation*}%
The order on $Y$ may be extended to $S$ in two ways. The first is
lexicographic. The second order is length-lexicographic, where first we
compare lengths of two sequences and then order lexicographically those
sequences of the same length.

Let $\mathbf{k}$ be field and $V/\mathbf{k}$ be the vector space with basis $%
S$. If $u\in S$ then let $uV$ denote the linear span of $uS$ and $V_{n}$
denote the $\mathbf{k}$-span of $S_{n}$ for all $n\geq 0$. Then, 
\begin{equation*}
V=\mathbf{k}\phi \oplus \dsum\limits_{y\in Y}yV=\cup _{n\geq 0}V_{n}\text{.}
\end{equation*}%
On repeating the first decomposition, we obtain subspaces $uV$ of $V$ where $%
yuV\leq uV$ for all $y\in Y,u\in S$; thus, the subspaces $uV$ form a
one-rooted tree $\mathbf{T}_{m}$. The space $V$ can also be seen as the
tensor space $\mathbf{k}\oplus \dsum\limits_{i\geq 0}\left( \otimes
^{i}U\right) $ where $U=V_{1}$, $V_{i+1}=\otimes ^{i}U$ for $i\geq 1$.

We consider endomorphisms $L$ of $V$ which leave $V_{n}$ invariant for all $%
n $. Thus $L$ corresponds to an element of direct product ring $E\left( 
\mathbf{T}_{m},\mathbf{k}\right) =\dprod_{m\geq 0}End_{\mathbf{k}}\left(
V_{n}\right) $ which is clearly residually finite dimensional. Another form
for $L$ which interests us is the following%
\begin{eqnarray*}
\left( \phi \right) L &=&L_{\phi }\phi , \\
\left( y.u\right) L &=&\dsum\limits_{y^{\prime }\in Y}y^{\prime }.\left(
u\right) L_{y,y^{\prime }}
\end{eqnarray*}%
where $L_{\phi }\in \mathbf{k}$ and where $L_{y,y^{\prime }}\in E\left( 
\mathbf{T}_{m},\mathbf{k}\right) $ for each pair $y,y^{\prime }\in Y$.

Using the length-lexicographic order on $S$ and writing the vectors in row
form, we get a matrix representation of $L\in E\left( \mathbf{T}_{m},\mathbf{%
k}\right) $ in block diagonal form 
\begin{equation*}
\lceil L\rceil =\left( 
\begin{array}{cccc}
L_{\phi } & 0 & 0 & . \\ 
0 & L_{m\times m} & 0 & . \\ 
0 & 0 & L_{m^{2}\times m^{2}} & . \\ 
. & . & . & .%
\end{array}%
\right)
\end{equation*}%
where $L_{m^{n}\times m^{n}}$ is the matrix of the transformation $L_{n}$
induced by $L$ on $V_{n}$ for all $n$. The ring of matrices $\lceil L\rceil $
is denoted by $M\left( \mathbf{T}_{m},\mathbf{k}\right) $.

On the other hand, on using the lexicographically ordered basis $S$ we
obtain the matrix representation of $L$ as an $\left( m+1\right) \times
\left( m+1\right) $ block form

\begin{equation*}
\left[ L\right] =\left( 
\begin{array}{cc}
L_{\phi } & 0 \\ 
0 & \left( \left[ L_{y_{i},y_{j}}\right] \right) _{m\times m}%
\end{array}%
\right)
\end{equation*}
The set of matrices $\left[ L\right] $ is denoted by $\mathbf{M}\left( 
\mathbf{T}_{m},\mathbf{k}\right) $.

We note that in the first matrix representation $\lceil L\rceil $ the block $%
L_{m\times m}$ has for its entries the values of $L_{y,y^{\prime }}$ at $%
\phi $ and the entries of the block $L_{m^{2}\times m^{2}}$ the values of $%
\left( L_{y,y^{\prime }}\right) _{y^{\prime \prime },y^{\prime \prime \prime
}}$ at $\phi $ and so on.

Define the set of states of $L$ as 
\begin{equation*}
Q(L)=\left\{ 
\begin{array}{c}
\left[ L\right] , \\ 
L_{\phi },\left[ L_{y,y^{\prime }}\right] , \\ 
\left( L_{y}\right) _{\phi },\left[ \left( L_{y,y^{\prime }}\right)
_{y^{\prime \prime },y^{\prime \prime \prime }}\right] ,...\text{.}%
\end{array}%
\right\} \text{.}
\end{equation*}%
\textbf{Example 1}. Automorphisms of the binary tree

Write $y_{1}=\dot{0},y_{2}=\dot{1}$. Then automorphism of the tree $\sigma
:\phi \rightarrow \phi ,\dot{0}u\leftrightarrow \dot{1}u$ is and in $M\left( 
\mathbf{T}_{2},\mathbf{k}\right) $ as%
\begin{equation*}
\lceil \sigma \rceil =\left( 
\begin{array}{cccc}
1 & 0 & 0 & . \\ 
0 & \sigma _{2\times 2} & 0 & . \\ 
0 & 0 & \sigma _{4\times 4} & . \\ 
. & . & . & .%
\end{array}%
\right)
\end{equation*}%
where $\sigma _{2\times 2}=\left( 
\begin{array}{cc}
0 & 1 \\ 
1 & 0%
\end{array}%
\right) ,\sigma _{2^{i+1}\times 2^{i+1}}=\sigma _{2^{i}\times 2^{i}}\otimes
I_{2\times 2}$ for $i\geq 1$. Also, $\sigma $ is represented in $\mathbf{M}%
\left( \mathbf{T}_{2},\mathbf{k}\right) $ as 
\begin{equation*}
\left[ \sigma \right] =\left( 
\begin{array}{ccc}
1 & 0 & 0 \\ 
0 & 0 & I \\ 
0 & I & 0%
\end{array}%
\right)
\end{equation*}

A general automorphism of the binary tree has the form $\alpha =\left(
\alpha _{0},\alpha _{1}\right) $ (inactive) or $\alpha =\left( \alpha
_{0},\alpha _{1}\right) \sigma $ (active) in $M\left( \mathbf{T}_{2},\mathbf{%
k}\right) $ as 
\begin{equation*}
\lceil \alpha \rceil =\left( 
\begin{array}{cccc}
1 & 0 & 0 & . \\ 
0 & I_{2} & 0 & . \\ 
0 & 0 & \alpha _{4\times 4} & . \\ 
. & . & . & .%
\end{array}%
\right) \text{ or }\left( 
\begin{array}{cccc}
1 & 0 & 0 & . \\ 
0 & \sigma _{2\times 2} & 0 & . \\ 
0 & 0 & \alpha _{4\times 4} & . \\ 
. & . & . & .%
\end{array}%
\right)
\end{equation*}%
where in the first case $\alpha _{4\times 4}=\left( 
\begin{array}{cc}
\beta _{2\times 2} & 0 \\ 
0 & \gamma _{2\times 2}%
\end{array}%
\right) $ and in the second case $\alpha _{4\times 4}=\left( 
\begin{array}{cc}
\beta _{2\times 2} & 0 \\ 
0 & \gamma _{2\times 2}%
\end{array}%
\right) \left( \sigma _{2\times 2}\otimes I_{2}\right) $. In $\mathbf{M}%
\left( \mathbf{T}_{2},\mathbf{k}\right) $, $\alpha $ is represented as 
\begin{equation*}
\left[ \alpha \right] =\left( 
\begin{array}{ccc}
1 & 0 & 0 \\ 
0 & \left[ \alpha _{0}\right] & 0 \\ 
0 & 0 & \left[ \alpha _{1}\right]%
\end{array}%
\right) \text{ or }\left( 
\begin{array}{ccc}
1 & 0 & 0 \\ 
0 & 0 & \left[ \alpha _{0}\right] \\ 
0 & \left[ \alpha _{1}\right] & 0%
\end{array}%
\right)
\end{equation*}%
according to the two cases above.

We extend the definitions of \textit{finite-state} and of \textit{%
functionally recursive }automorphisms of rooted regular trees in \cite%
{brun-sid97} to $\mathbf{M}\left( \mathbf{T}_{m},\mathbf{k}\right) $.

\begin{definition}
1. A matrix $L\in \mathbf{M}\left( \mathbf{T}_{m},\mathbf{k}\right) $ is
said to be \textit{finite-state }provided $Q\left( L\right) $ is finite.%
\newline
2. A finite set of matrices $\left\{ N_{i}\mid i=1...,v\right\} $ in $%
\mathbf{M}\left( \mathbf{T}_{m},\mathbf{k}\right) $ is said to be \textit{%
functionally recursive} provided each of the states $\left( N_{i}\right)
_{y,y^{\prime }}$ is the value of some polynomial in non-commuting variables
with constants from $M\left( m+1,\mathbf{k}\right) $, evaluated at the $%
N_{j} $'s. A matrix $L$ is said to \textit{functionally recursive provided
it is an element of some functionally recursive set. }
\end{definition}

A finite-state matrix is functionally recursive. Both sets of matrices form $%
\mathbf{k}$-algebras. The notion of finite-state transformation may be
extended by assuming the space generated by the states to be finite
dimensional.

Our constructions in the large algebra $M\left( \mathbf{T}_{m},\mathbf{k}%
\right) $ often deal with sequences of transformations which are convergent
in some sense, or have a special growth type. We formalize a notion of
convergence in a more general setting. Let $\left\{ W_{i}/\mathbf{k}\mid
i\geq 1\right\} $ be a sequence of finite dimensional vector spaces with $%
dim(W_{i})=m_{i}$, let $L_{i}\in End_{\mathbf{k}}\left( W_{i}\right) $ and
define $r_{i}=\frac{rank\left( L_{i}\right) }{m_{i}}$. Under the usual
operations, the set of sequences $\left\{ L_{i}\mid i\geq 1\right\} $ form a 
$\mathbf{k}$-algebra. We say that the sequence $\left\{ L_{i}\mid i\geq
1\right\} $ \textit{converges} provided the sequence $r=\left\{ r_{i}\mid
i\geq 1\right\} $ converges. It is direct to show that the set of convergent
sequences of transformations is a $\mathbf{k}$-algebra. Clearly, if the $%
L_{i}$%
%TCIMACRO{\U{b4}}%
%BeginExpansion
\'{}%
%EndExpansion
s are invertible then $r$ is the constant sequence $1$. Also, those
sequences of transformations whose $r$ converges to $0$ form a $\mathbf{k}$%
-algebra.

One notion of growth for an automorphism $\alpha $ of a tree counts the
number of active states in $\left\{ \alpha _{u}\mid u\in S_{n}\right\} $ for
each $n$ \cite{sid00}. We give here an additive form of this version of
growth. Consider $V$ as the tensor space $\mathbf{k}\oplus
\dsum\limits_{i\geq 0}\left( \otimes ^{i}U\right) $, where $V_{0}=\mathbf{k}$
and $V_{i+1}=\otimes ^{i}U$ for $i\geq 1$. Let $id\left( U\right) $ be the
identity transformation on $U$. Then $End_{\mathbf{k}}\left( V_{i}\right) $
embeds in $End_{\mathbf{k}}\left( V_{i+1}\right) $ via the map $%
L_{i}\rightarrow L_{i}\otimes id\left( U\right) $. Define $\partial
L_{i+1}=L_{i+1}-L_{i}\otimes id\left( U\right) $ for $i\geq 0$. We define
the \textit{growth function} associated to the sequence $\left\{ \partial
L_{i}\mathbf{\mid i\geq }1\right\} $ as $f:i\rightarrow rank\left( \partial
L_{i}\right) $. The set of sequences $\left\{ L_{i}\mathbf{\mid i\geq }%
1\right\} $ of a given growth type form a $\mathbf{k}$-algebra.

\section{A free $2$-generated ring}

Let $s=\left( 
\begin{array}{ccc}
1 & 0 & 0 \\ 
0 & I & 0 \\ 
0 & 0 & 2s%
\end{array}%
\right) ,t=\left( 
\begin{array}{ccc}
0 & 0 & 0 \\ 
0 & 0 & 2s \\ 
0 & 0 & 2t%
\end{array}%
\right) \in \mathbf{M}\left( \mathbf{T}_{2},\mathbb{Q}\right) $. Write 
\begin{equation*}
s=\left( 
\begin{array}{cc}
1 & 0 \\ 
0 & 2s%
\end{array}%
\right) ,t=\left( 
\begin{array}{cc}
0 & 2s \\ 
0 & 2t%
\end{array}%
\right)
\end{equation*}%
for short. Then $s,t$ are functionally recursive and have exponential
growth. Let $R$ be the ring generated by $s$ and $t$.

Given a monomial $m=s^{i_{1}}t^{j_{1}}...s^{i_{l}}t^{j_{l}}$, where $%
i_{1},...,i_{l},j_{1},...,j_{l}$ are non-negative integers, let $\left\vert
m\right\vert $ denote its formal length $i_{1}+j_{1}+...+i_{l}+j_{l}$. We
find that for $i_{1}\geq 0,j_{1}\geq 1$ and $m^{\prime
}=s^{i_{2}}t^{j_{2}}...s^{i_{l}}t^{j_{l}}$,%
\begin{equation*}
m=\left( 
\begin{array}{cc}
0 & 2^{\left\vert m\right\vert -i_{1}}st^{j_{1}-1}m^{\prime } \\ 
0 & 2^{\left\vert m\right\vert }m%
\end{array}%
\right) \text{.}
\end{equation*}

\begin{lemma}
The generator $s$ is transcendental in $R$.
\end{lemma}

\begin{proof}
Suppose not. Let $p\left( x\right) =\sum_{0\leq i\leq l}a_{i}x^{i}$ be a
non-zero polynomial of minimal degree $l$ such that $p\left( s\right) =0$.
Note that $s$ is invertible in $M\left( \mathbb{Q},\mathbf{T}_{2}\right) $
and so $a_{0}\not=0$. Now, from the matrix form%
\begin{equation*}
p\left( s\right) =\left( 
\begin{array}{cc}
\sum_{0\leq i\leq l}a_{i} & 0 \\ 
0 & \sum_{0\leq i\leq l}2^{i}a_{i}s^{i}%
\end{array}%
\right)
\end{equation*}%
we have 
\begin{equation*}
\sum_{0\leq i\leq l}a_{i}=0,\sum_{0\leq i\leq l}2^{i}a_{i}s^{i}=0\text{.}
\end{equation*}%
Therefore, 
\begin{equation*}
\sum_{0\leq i\leq l}2^{i}a_{i}s^{i}-\sum_{0\leq i\leq
l}a_{i}s^{i}=\sum_{1\leq i\leq l}\left( 2^{i}-1\right) a_{i}s^{i}=0\text{.}
\end{equation*}%
It follows that $\left( 2^{i}-1\right) a_{i}=0$ and so, $a_{i}=0$ for all $%
1\leq i\leq l.$ Thus, $\sum_{0\leq i\leq l}a_{i}=a_{0}=0$; a contradiction.
\end{proof}

A word $w$ in $R$ has the form%
\begin{eqnarray*}
w &=&\sum_{0\leq i\leq l}a_{i}s^{i}+\sum_{1\leq k\leq n}b_{k}m_{k}\text{,} \\
m_{k} &=&s^{i_{1}\left( m_{k}\right) }t^{j_{1}\left( m_{k}\right)
}m_{k}^{\prime },\text{ where }i_{1}\left( m_{k}\right) \geq 0,j_{1}\left(
m_{k}\right) >0\text{.}
\end{eqnarray*}%
The matrix form of $w$ is 
\begin{equation*}
w=\left( 
\begin{array}{cc}
\sum_{0\leq i\leq l}a_{i} & \sum_{1\leq k\leq n}2^{\left\vert
m_{k}\right\vert -i_{1}\left( m_{k}\right) }b_{k}st^{j_{1}\left(
m_{k}\right) -1}m_{k}^{\prime } \\ 
0 & \sum_{0\leq i\leq l}2^{i}a_{i}s^{i}+\sum_{1\leq k\leq n}2^{\left\vert
m_{k}\right\vert }b_{k}m_{k}%
\end{array}%
\right) \text{.}
\end{equation*}%
\textbf{Proof of Theorem 1.}

Suppose there exists a non-trivial word $w$ such that $w=0$ in $R$. By the
previous lemma, we may assume $b_{n}\not=0$; choose$\ w$ having minimal $v$
such that $b_{n}\not=0$. Then,%
\begin{eqnarray*}
\sum_{0\leq i\leq l}a_{i} &=&0, \\
\sum_{0\leq i\leq l}2^{i}a_{i}s^{i}+\sum_{1\leq k\leq n}2^{\left\vert
m_{k}\right\vert }b_{k}m_{k} &=&0, \\
\sum_{1\leq k\leq n}2^{\left\vert m_{k}\right\vert -i_{1}\left( m_{k}\right)
}b_{k}st^{j_{1}\left( m_{k}\right) -1}m_{k}^{\prime } &=&0.
\end{eqnarray*}%
We eliminate the leading term of $w$:%
\begin{equation*}
w_{22}-2^{\left\vert m_{n}\right\vert }w=\sum \left( 2^{i}-2^{\left\vert
m_{n}\right\vert }\right) a_{i}s^{i}+\sum_{1\leq k\leq l-1}\left(
2^{\left\vert m_{k}\right\vert }-2^{\left\vert m_{n}\right\vert }\right)
b_{k}m_{k}\text{.}
\end{equation*}%
By the minimality of $w$, 
\begin{equation*}
\left( 2^{i}-2^{\left\vert m_{n}\right\vert }\right) a_{i}=0,\left(
2^{\left\vert m_{k}\right\vert }-2^{\left\vert m_{n}\right\vert }\right)
b_{k}=0\text{.}
\end{equation*}%
If $a_{i}\not=0$ for some $i$, then $i=\left\vert m_{n}\right\vert $ and $%
a_{j}=0$ for all $j\not=i$; but, as $\sum_{0\leq j\leq u}a_{j}=0$, we have a
contradiction. Hence, $w=\sum_{1\leq k\leq n}b_{k}m_{k}$ and $\left\vert
m_{k}\right\vert =\left\vert m_{n}\right\vert $ for all $k$, and as $s$ is
invertible, 
\begin{equation*}
\sum_{1\leq k\leq n}2^{\left\vert m_{n}\right\vert -i_{1}\left( m_{k}\right)
}b_{k}t^{j_{1}\left( m_{k}\right) -1}m_{k}^{\prime }=0\text{.}
\end{equation*}%
Suppose there exist $k\not=l$ such that $t^{j_{1}\left( m_{k}\right)
-1}m_{k}^{\prime }=t^{j_{1}\left( m_{l}\right) -1}m_{l}^{\prime }$. It
follows that, $t^{j_{1}\left( m_{k}\right) }m_{k}^{\prime }=t^{j_{1}\left(
m_{l}\right) }m_{l}^{\prime }$, and 
\begin{equation*}
j_{1}\left( m_{k}\right) +\left\vert m_{k}^{\prime }\right\vert =j_{1}\left(
m_{l}\right) +\left\vert m_{l}^{\prime }\right\vert \text{.}
\end{equation*}%
Now, since $m_{k}=s^{i_{1}\left( m_{k}\right) }t^{j_{1}\left( m_{k}\right)
}m_{k}^{\prime }$ and $\left\vert m_{k}\right\vert =\left\vert
m_{l}\right\vert $ we have%
\begin{equation*}
\left\vert m_{k}\right\vert =i_{1}\left( m_{k}\right) +j_{1}\left(
m_{k}\right) +\left\vert m_{k}^{\prime }\right\vert ,
\end{equation*}%
\begin{eqnarray*}
i_{1}\left( m_{k}\right) +j_{1}\left( m_{k}\right) +\left\vert m_{k}^{\prime
}\right\vert &=&i_{1}\left( m_{l}\right) +j_{1}\left( m_{l}\right)
+\left\vert m_{l}^{\prime }\right\vert , \\
i_{1}\left( m_{k}\right) &=&i_{1}\left( m_{l}\right) , \\
t^{j_{1}\left( m_{k}\right) }s^{i_{1}\left( m_{k}\right) }m_{k}^{\prime }
&=&t^{j_{1}\left( m_{l}\right) }s^{i_{1}\left( m_{l}\right) }m_{l}^{\prime },
\\
m_{k} &=&m_{l}\text{.}
\end{eqnarray*}%
Therefore, $k=l$; a contradiction. Hence, by minimality of $w$, we obtain $%
2^{\left\vert m_{n}\right\vert -i_{1}\left( m_{k}\right) }b_{k}=0$ and $%
b_{k}=0$ for all $k$; a final contradiction.

\section{A $2$-generated monomial ring}

Let $s=\left( 
\begin{array}{ccc}
0 & 0 & 0 \\ 
0 & 0 & 0 \\ 
0 & I & 0%
\end{array}%
\right) ,t=\left( 
\begin{array}{ccc}
0 & 0 & 0 \\ 
0 & 0 & t \\ 
0 & 0 & s%
\end{array}%
\right) \in \mathbf{M}\left( \mathbf{T}_{2},\mathbb{Q}\right) $. Write%
\begin{equation*}
s=\left( 
\begin{array}{cc}
0 & 0 \\ 
1 & 0%
\end{array}%
\right) ,t=\left( 
\begin{array}{cc}
0 & t \\ 
0 & s%
\end{array}%
\right)
\end{equation*}
for short. Then $s,t$ are finite-state transformations and both have bounded
growth. Let $R$ be the ring generated by $s$ and $t$.

Before stating results about the ring $R$, we introduce some notation. For a
general ring $A$ generated by a set $X$, let $\mu \left( A\right) $ denote
the multiplicative semigroup generated by $X$.

Given $u\in \mu \left( A\right) $ we let $|u|$ denote the formal length of $%
u $ as a word in $X$. For $w\in A$, let $\delta \left( w\right) $ denote the
nilpotency degree of $w$; if $w$ is not nilpotent, then let $\delta \left(
w\right) $ be infinite.

Let $ran\left( A\right) ,lan\left( A\right) $ the the right and left
annihilators of $A$, respectively. Also, define the ascending left
annihilator series $L_{1}\left( A\right) =lan\left( A\right) ,L_{i+1}\left(
A\right) =lan\left( A,L_{i}\left( A\right) \right) $ (that is, $\frac{%
L_{i+1}\left( A\right) }{L_{i}\left( A\right) }=lan\left( \frac{A}{%
L_{i}\left( A\right) }\right) $), $L_{\infty }\left( A\right)
=\dbigcup\limits_{i\geq 1}L_{i}\left( A\right) $ and $\overline{R}=\frac{R}{%
L_{\infty }\left( R\right) }$.

We will prove in the next sections various properties of $R$. Some of these
are gathered below.

\begin{theorem}
(i) The ring $R$ is torsion-free and has infinite $\mathbb{Z}$-rank.\newline
(ii) The element $w=s+t$ is transcendental over $\mathbb{Z}$.$\newline
$(iii) The ideal $ran\left( R\right) =0$, and the ideals $L_{i}$ $\left(
R\right) $ are nilpotent$\newline
$ with $L_{i+1}\left( R\right) >L_{i}$ $\left( R\right) $ for all $i\geq 1$.%
\newline
(iv) The ideals $ran\left( \overline{R}\right) =lan\left( \overline{R}%
\right) =0$. Furthermore, $\overline{R}$ has infinite $\mathbb{Z}$-rank $%
\newline
$(v) The semigroup $\mu \left( R\right) $ is nil of degree $5$ and $\mu
\left( \overline{R}\right) $ is nil of degree $4$.\newline
(vi) The ring $R$ modulo the ideal generated by $\left( tst\right) ^{3}$ has
finite $\mathbb{Z}$-rank.\newline
(vii) The ring $R$ affords the monomial presentation in Theorem 2.
\end{theorem}

\subsection{First structural information}

Let $\widehat{R}=\mathbb{Z}1\oplus R$. Then, $R\leq M_{2\times 2}\left( 
\widehat{R}\right) $ and we have the following decompositions 
\begin{equation*}
R=\widehat{R}s\oplus \widehat{R}t=s\widehat{R}\oplus t\widehat{R}\text{.}
\end{equation*}%
Let $T=$ $\widehat{R}t=\left\langle t,st\right\rangle $. Then $T$ is a
proper subring of $R$ and the projection $\varphi :t=\left( 
\begin{array}{cc}
0 & t \\ 
0 & s%
\end{array}%
\right) \rightarrow s,st=\left( 
\begin{array}{cc}
0 & 0 \\ 
0 & t%
\end{array}%
\right) \rightarrow t$ extends to an epimorphism $\varphi :T\rightarrow R$;
in other words, the ring $R$ is recurrent. Therefore, we conclude from $R=%
\widehat{R}s\oplus T$ that $R$ has infinite $\mathbb{Z}$-rank. We decompose $%
R$ in terms of $T$.

\begin{lemma}
$R=\mathbb{Z}s\oplus st\widehat{T}s\oplus st.\widehat{T}\oplus t.\widehat{T}%
.s\oplus t.\widehat{T}$.

\begin{proof}
\begin{eqnarray*}
R &=&s\widehat{R}\oplus t\widehat{R}=s\left( \mathbb{Z}\oplus R\right)
\oplus t\left( \mathbb{Z}\oplus R\right) \\
&=&\mathbb{Z}s\oplus s\left( \widehat{R}s\oplus \widehat{R}t\right) \oplus 
\mathbb{Z}t\oplus t\left( \widehat{R}s\oplus \widehat{R}t\right) \\
&=&\mathbb{Z}s\oplus s\widehat{R}s\oplus s\widehat{R}t\oplus \mathbb{Z}%
t\oplus t\widehat{R}s\oplus t\widehat{R}t \\
&=&\mathbb{Z}s\oplus s\left( \mathbb{Z}\oplus \widehat{R}s\oplus \widehat{R}%
t\right) s\oplus s\left( \mathbb{Z}\oplus s\widehat{R}\oplus t\widehat{R}%
\right) t\oplus \mathbb{Z}t\oplus \\
&&t\left( \mathbb{Z}\oplus \widehat{R}s\oplus \widehat{R}t\right) s\oplus t%
\widehat{R}t \\
&=&\mathbb{Z}s\oplus s\widehat{R}ts\oplus \left( \mathbb{Z}st\oplus st%
\widehat{R}t\right) \oplus \mathbb{Z}t\oplus \left( \mathbb{Z}ts\oplus t%
\widehat{R}ts\right) \oplus t\widehat{R}t \\
&=&\mathbb{Z}s\oplus s\left( \mathbb{Z}\oplus s\widehat{R}\oplus t\widehat{R}%
\right) ts\oplus \oplus st.\widehat{T}\oplus t.\widehat{T}.s\oplus t.%
\widehat{T} \\
&=&\mathbb{Z}s\oplus st\widehat{T}s\oplus st.\widehat{T}\oplus t.\widehat{T}%
.s\oplus t.\widehat{T}\text{.}
\end{eqnarray*}
\end{proof}
\end{lemma}

The above decomposition provides us with five monomial types.

\begin{lemma}
Elements of $\mu \left( R\right) $ have the following five matrix types and
nilpotency degrees:%
\begin{eqnarray*}
s &=&\left( 
\begin{array}{cc}
0 & 0 \\ 
1 & 0%
\end{array}%
\right) ,\delta \left( s\right) =2; \\
t.v\left( t,st\right) &=&\left( 
\begin{array}{cc}
0 & t.v(s,t) \\ 
0 & s.v\left( s,t\right)%
\end{array}%
\right) ,\delta \left( t.v\left( t,st\right) \right) =\delta \left(
s.v\left( s,t\right) \right) \text{ or }\delta \left( s.v\left( s,t\right)
\right) +1; \\
st.v\left( t,st\right) &=&\left( 
\begin{array}{cc}
0 & 0 \\ 
0 & t.v\left( s,t\right)%
\end{array}%
\right) ,\delta \left( st.v\left( t,st\right) \right) =\delta \left(
t.v\left( s,t\right) \right) ; \\
t.v\left( t,st\right) .s &=&\left( 
\begin{array}{cc}
t.v(s,t) & 0 \\ 
s.v\left( s,t\right) & 0%
\end{array}%
\right) ,\delta \left( t.v\left( t,st\right) .s\right) =\delta \left(
s.v\left( s,t\right) \right) \text{ or }\delta \left( s.v\left( s,t\right)
\right) +1; \\
st.w\left( t,st\right) .s &=&\left( 
\begin{array}{cc}
0 & 0 \\ 
t.v\left( s,t\right) & 0%
\end{array}%
\right) ,\delta \left( st.v\left( t,st\right) .s\right) =\delta \left(
t.v\left( s,t\right) \right) \text{.}
\end{eqnarray*}
\end{lemma}

The following short list of monomials and their nilpotency degrees will be
helpful in later computations.

\begin{lemma}
The monomials of length at most $3$ are%
\begin{eqnarray*}
s &=&\left( 
\begin{array}{cc}
0 & 0 \\ 
1 & 0%
\end{array}%
\right) ,t=\left( 
\begin{array}{cc}
0 & t \\ 
0 & s%
\end{array}%
\right) ,s^{2}=0,ts=\left( 
\begin{array}{cc}
t & 0 \\ 
s & 0%
\end{array}%
\right) , \\
st &=&\left( 
\begin{array}{cc}
0 & 0 \\ 
0 & t%
\end{array}%
\right) ,t^{2}=\left( 
\begin{array}{cc}
0 & ts \\ 
0 & 0%
\end{array}%
\right) ,tst=\left( 
\begin{array}{cc}
0 & t^{2} \\ 
0 & st%
\end{array}%
\right) , \\
sts &=&\left( 
\begin{array}{cc}
0 & 0 \\ 
t & 0%
\end{array}%
\right) ,t^{2}s=\left( 
\begin{array}{cc}
ts & 0 \\ 
0 & 0%
\end{array}%
\right) ,st^{2}=\left( 
\begin{array}{cc}
0 & 0 \\ 
0 & ts%
\end{array}%
\right) ,t^{3}=0\text{.}
\end{eqnarray*}%
The nilpotency degrees of these elements are as follows:%
\begin{eqnarray*}
\delta \left( s\right) &=&2,\delta \left( t\right) =3,\delta \left(
ts\right) =4, \\
\delta \left( st\right) &=&3,\delta \left( t^{2}\right) =2,\delta \left(
tst\right) =4, \\
\delta \left( sts\right) &=&2,\delta \left( t^{2}s\right) =4,\delta \left(
st^{2}\right) =4\text{.}
\end{eqnarray*}
\end{lemma}

We are now able to prove a number of important properties of $R$.

\begin{proposition}
(i) (Contraction property) Let $v\in \mu \left( R\right) $ and let $%
\left\vert v\right\vert \geq 3$. Then $v=\left( 
\begin{array}{cc}
0 & v_{12} \\ 
0 & v_{22}%
\end{array}%
\right) $ or $\left( 
\begin{array}{cc}
v_{11} & 0 \\ 
v_{21} & 0%
\end{array}%
\right) $ where $v_{ij}\in \mu \left( R\right) $ and $|v_{ij}|<|v|$.$\newline
$(ii) Let $u,v$ be non-zero monomials such that $u=v$ in $R$. Then, $u$ is
freely equal to $v$.$\newline
$(iii) $R$ is a monomial ring.$\newline
$(iv) Monomials of $R$ are nilpotent.$\newline
$(v) Let $u=s+t$. Then $u$ is transcendental over $\mathbb{Z}$.
\end{proposition}

\begin{proof}
(i) First check that the contraction property holds for monomials of length $%
3$. Now the assertion follows from Lemma 3, by induction on the length $|v|$.

(ii) Proceed by induction on max$\left\{ \left\vert u\right\vert ,\left\vert
v\right\vert \right\} $ and assume $\left\vert u\right\vert \leq \left\vert
v\right\vert $. From the matrix representations of $u,v$, we may assume $%
u=\left( 
\begin{array}{cc}
0 & u_{12} \\ 
0 & u_{22}%
\end{array}%
\right) ,v=\left( 
\begin{array}{cc}
0 & v_{12} \\ 
0 & v_{22}%
\end{array}%
\right) $. By the contraction property, we may assume furthermore, $%
\left\vert v\right\vert \leq 2$. The assertion follows by checking the list.

(iii) Let $v_{i}\not=0$ be formally distinct monomials and $n_{i}\in \mathbb{%
Z}$ $\left( i=1,...,k\right) $. If $\sum_{1\leq i\leq k}n_{i}v_{i}=0$, then
we will show that $n_{i}=0$ for all $i$.

Proceed by contradiction. Choose $k$ minimum and let $v_{i}\in \mu \left(
R\right) $ be such that $\sum_{k}|v_{i}|$ is minimum. We represent $v_{i}$
as in Lemma 3 above. Let

\begin{equation*}
I=\left\{ 1,...,k\right\} ,I_{1}=\left\{ i\mid \left( v_{i}\right)
_{12}\not=0\right\} ,I_{2}=\left\{ i\mid \left( v_{i}\right) _{12}=0\right\} 
\text{.}
\end{equation*}

Then, $\sum_{k}n_{i}v_{i}=\sum \left\{ n_{i}v_{i}\mid i\in I_{1}\right\}
+\sum \left\{ n_{i}v_{i}\mid i\in I_{2}\right\} =0$ implies 
\begin{equation*}
\sum \left\{ n_{i}\left( v_{i}\right) _{12}\mid i\in I_{1}\right\} =0\text{.}
\end{equation*}%
Suppose $I_{1}$ is nonempty. We have for all $i\in I_{1}$, $%
v_{i}=tw_{i}\left( t,st\right) =\left( 
\begin{array}{cc}
0 & t.w_{i}(s,t) \\ 
0 & s.w_{i}\left( s,t\right)%
\end{array}%
\right) $ where $\left\vert \left( v_{i}\right) _{12}\right\vert =\left\vert
t.w_{i}(s,t)\right\vert <\left\vert v_{i}\right\vert $ By induction on
length, for every $i$ in $I_{1}$, there exists $j\not=i$ in $I_{1}$ such
that $t.w_{i}(s,t)=t.w_{j}(s,t)$. Therefore, $t.w_{i}(s,t),t.w_{j}(s,t)$ are
freely equal and so, $w_{i}(s,t)=w_{j}(s,t)$ and $v_{i}=v_{j}$; a
contradiction.

Similarly, it follows that $\left\{ i\in I\mid \left( v_{i}\right)
_{11}\not=0\right\} $ is empty. Hence, for all $i$ in $I$ we have $%
v_{i}=st.w_{i}\left( t,st\right) .s=\left( 
\begin{array}{cc}
0 & 0 \\ 
t.w_{i}\left( s,t\right) & 0%
\end{array}%
\right) $ or $v_{i}=st.w_{i}\left( t,st\right) =\left( 
\begin{array}{cc}
0 & 0 \\ 
0 & t.w_{i}\left( s,t\right)%
\end{array}%
\right) $.

We repeat the argument above for the partition%
\begin{equation*}
J_{1}=\left\{ i\in I\mid \left( v_{i}\right) _{21}\not=0\right\}
,J_{2}=\left\{ i\in I\mid \left( v_{i}\right) _{22}=0\right\} \text{.}
\end{equation*}

(iv) Let $v\in \mu \left( R\right) $. Then, $v^{4}=0$ for $|v|\leq 2$.
Suppose $3\leq \left\vert v\right\vert $ and $v=\left( 
\begin{array}{cc}
0 & v_{12} \\ 
0 & v_{22}%
\end{array}%
\right) $. Then, as $\left\vert v_{22}\right\vert <|v|$ we have $v_{22}$ is
nilpotent, say $\delta \left( v_{22}\right) =n$. Now since $v^{i}=\left( 
\begin{array}{cc}
0 & v_{12}v_{22}^{i-1} \\ 
0 & v_{22}^{i}%
\end{array}%
\right) $ for all $i\geq 1$, we conclude $v^{n+1}=0$. The case where $%
v=\left( 
\begin{array}{cc}
v_{11} & 0 \\ 
v_{21} & 0%
\end{array}%
\right) $ is similar.

(v) Let $u=s+t$. For all $n\geq 1$, $u^{n}$ is the sum of all non-zero
monomials of length $n$. Since $\mu \left( R\right) $ is infinite and $R$ is
a monomial ring, it follows that $u$ is transcendental.
\end{proof}

Since $R$ is a monomial ring, it follows clearly that $R$ is torsion-free.

\subsection{Right and Left Annihilators}

A subring $U$ of $R$ is called a \textit{monomially closed} provided if $%
w\in U$ and $w=\sum u_{i}$ where $u_{i}\in \mu \left( R\right) $ then $%
u_{i}\in \mu \left( U\right) $. If $u\in R$ then let $ran\left( u\right) $
and $lan\left( u\right) $ denote the right and left annihilator of $u$,
respectively.

\begin{lemma}
Let $u$ be a homogenous element in $R$. Then, $ran\left( u\right) $ and $%
lan\left( u\right) $ are monomially closed.
\end{lemma}

\begin{proof}
The element $u=\sum n_{i}u_{i}$ where the $u_{i}$%
%TCIMACRO{\U{b4}}%
%BeginExpansion
\'{}%
%EndExpansion
s are distinct monomials of the same length. Suppose $w=\sum m_{j}w_{j}\in
ran\left( u\right) $ where the $w_{i}$ are monomials and the sum is reduced.
Then, $uw=\sum n_{i}m_{j}u_{i}w_{j}=0$ and so, $u_{i}w_{j}$ are distinct or
null.
\end{proof}

\begin{lemma}
(i) $lan\left( s\right) =\widehat{R}s$, (ii) $ran\left( R\right) =\left\{
0\right\} $.
\end{lemma}

\begin{proof}
Since $lan\left( s\right) $ and $ran\left( R\right) $ are monomially closed,
it suffices to describe the monomials in these subrings. We will use in our
analysis the five types of monomials $u$ in Lemma 3.

(i) We note that $u\in lan\left( s\right) \cap \mu \left( R\right) $ if and
only if $u=s,t.w\left( t,st\right) .s,st.w\left( t,st\right) .s$ and these
clearly generate $\widehat{R}s$.

(ii) Suppose $s.u=t.u=0$. We check that $u=0$ or $\left\vert u\right\vert
\geq 3$. Suppose $u=t.w\left( t,st\right) =\left( 
\begin{array}{cc}
0 & tw(s,t) \\ 
0 & sw\left( s,t\right)%
\end{array}%
\right) $. Then, $t.w(s,t)=0,ts.w\left( s,t\right) =0$. It follows that $%
s.w\left( s,t\right) \in ran\left( R\right) $ and by induction, $s.w\left(
s,t\right) =0$. Therefore, we have $t.w(s,t)=0,s.w\left( s,t\right) =0$ and $%
w(s,t)\in ran\left( R\right) $, $w(s,t)=0$.
\end{proof}

\begin{proposition}
The ideals in $\left\{ L_{i}\left( R\right) \mid i\geq 1\right\} $ are
nilpotent and distinct. The quotient ring $\frac{R}{L_{\infty }\left(
R\right) }$ has infinite $\mathbb{Z}$-rank.
\end{proposition}

\begin{proof}
Note that $u\in \mu \left( L_{i}\right) $ if and only if $uw=0$ for all $%
w\in \mu \left( R\right) $ of length $i$. This shows that $L\left( R\right)
^{2}=\left\{ 0\right\} $ and more generally, $L_{i}\left( R\right)
^{1+i}=\left\{ 0\right\} $.

Define the sequence 
\begin{equation*}
x_{1}=ts,\text{ }x_{i+1}\left( s,t\right) =x_{i}\left( s,t\right) ^{\varphi
}=x_{i}\left( t,st\right) \text{ for all }i\geq 1
\end{equation*}%
and

\begin{equation*}
y_{i}=x_{i}^{3}\text{ for all }i\geq 1\text{.}
\end{equation*}%
Then 
\begin{eqnarray*}
x_{2} &=&st^{2}=\left( 
\begin{array}{cc}
0 & 0 \\ 
0 & ts%
\end{array}%
\right) ,x_{3}=t\left( st\right) ^{2}=\left( 
\begin{array}{cc}
0 & 0 \\ 
0 & st^{2}%
\end{array}%
\right) ,x_{i+1}=\left( 
\begin{array}{cc}
0 & 0 \\ 
0 & x_{i}%
\end{array}%
\right) \text{ for all }i\geq 1\text{;} \\
y_{1} &=&\left( ts\right) ^{3}=\left( 
\begin{array}{cc}
0 & 0 \\ 
st^{2} & 0%
\end{array}%
\right) ,y_{2}=\left( 
\begin{array}{cc}
0 & 0 \\ 
0 & v_{1}%
\end{array}%
\right) ,y_{i+1}=\left( 
\begin{array}{cc}
0 & 0 \\ 
0 & y_{i}%
\end{array}%
\right) \text{ for all }i\geq 1\text{.}
\end{eqnarray*}

Since $y_{1}s=0$, $y_{1}t=0$, we conclude that $y_{1}\in L_{1}=lan\left(
R\right) $. Furthermore, $y_{2}s=\left( 
\begin{array}{cc}
0 & 0 \\ 
y_{1} & 0%
\end{array}%
\right) ,v_{2}t=\left( 
\begin{array}{cc}
0 & 0 \\ 
0 & y_{1}s%
\end{array}%
\right) $ and so, $y_{2}\in L_{2}\backslash L_{1}$.

More generally, $y_{i}\in L_{i}\backslash L_{i-1}$ $\left( i\geq 1\right) $
follows from 
\begin{eqnarray*}
\left( 
\begin{array}{cc}
L_{i} & 0 \\ 
0 & 0%
\end{array}%
\right) \cap R\oplus \left( 
\begin{array}{cc}
0 & 0 \\ 
L_{i} & 0%
\end{array}%
\right) \cap R &\leq &L_{i}, \\
\left( 
\begin{array}{cc}
0 & L_{i} \\ 
0 & 0%
\end{array}%
\right) \cap R\oplus \left( 
\begin{array}{cc}
0 & 0 \\ 
0 & L_{i}%
\end{array}%
\right) \cap R &\leq &L_{i+1}\text{.}
\end{eqnarray*}%
Suppose $\overline{R}=\frac{R}{L_{\infty }\left( R\right) }$ has finite
rank. Then, $\overline{R}$ is finite over $\mathbb{Z}_{2}$ and is therefore
nilpotent, say of degree $m$. Hence $\left( s+t\right) ^{m}$ $\in L_{\infty
}\left( R\right) $ and it follows that $s+t$ is nilpotent; a contradiction.
\end{proof}

\begin{lemma}
The right annihilators $ran\left( \frac{R}{L_{i}\left( R\right) }\right) $, $%
ran\left( \frac{R}{L_{\infty }\left( R\right) }\right) $ are null.
\end{lemma}

\begin{proof}
Let $w\in R$ such that $sw,tw\in lan\left( R\right) $. Then, $%
sws,tws,swt,twt $ are all zero. From $s\left( ws\right) =t\left( ws\right)
=0 $ we obtain $ws=0$ and likewise, $wt=0$. Therefore, $w\in lan\left(
R\right) $. Hence, $ran\left( \frac{R}{L_{1}\left( R\right) }\right)
=\left\{ 0\right\} $. By induction, $ran\left( \frac{R}{L_{i}\left( R\right) 
}\right) =\left\{ 0\right\} $ and $ran\left( \frac{R}{L_{\infty }\left(
R\right) }\right) =\left\{ 0\right\} $.
\end{proof}

\subsection{Powers of monomials}

The semigroup $\mu \left( R\right) $ is nil of degree at least $5$.

In order to calculate powers of monomials we need to evaluate words of the
form $\left( tw\right) \left( sw\right) ^{i}$ and $\left( sw\right) \left(
tw\right) ^{i}$, as in

\begin{eqnarray*}
\left( tw\left( t,st\right) \right) ^{k} &=&\left( 
\begin{array}{cc}
0 & tw \\ 
0 & sw%
\end{array}%
\right) ^{k}=\left( 
\begin{array}{cc}
0 & \left( tw\right) \left( sw\right) ^{k-1} \\ 
0 & \left( sw\right) ^{k}%
\end{array}%
\right) \text{,} \\
\left( tw\left( t,st\right) s\right) ^{k} &=&\left( 
\begin{array}{cc}
tw & 0 \\ 
sw & 0%
\end{array}%
\right) ^{k}=\left( 
\begin{array}{cc}
\left( tw\right) ^{k} & 0 \\ 
\left( sw\right) \left( tw\right) ^{k-1} & 0%
\end{array}%
\right) \text{.}
\end{eqnarray*}

\begin{proposition}
For any $w\in \mu \left( R\right) $ of length greater than $2$ and for any $%
k\geq 1$ the non-zero block entries of elements of 
\begin{equation*}
S(w,k)=\left\{ \left( tw\right) \left( sw\right) ^{k},\left( sw\right)
\left( tw\right) ^{k},\left( tw\right) \left( tsw\right) ^{k},\left(
tsw\right) \left( stw\right) ^{k},\left( tsw\right) \left( t^{2}w\right)
^{k}\right\}
\end{equation*}%
belong to $\cup \left\{ S(v,k)\mid \left\vert v\right\vert <\left\vert
w\right\vert \right\} $.
\end{proposition}

\begin{proof}
First, check that $S(w,k)=\left\{ 0\right\} $ for $w=sw^{\prime }$ or $%
t^{2}w^{\prime }$; thus, we need to verify the assertion only for $%
w=tv(t,st),tv(t,st)s$.

1. Let $w=tv(t,st)=\left( 
\begin{array}{cc}
0 & tv \\ 
0 & sv%
\end{array}%
\right) $. Then,%
\begin{eqnarray*}
sw &=&\left( 
\begin{array}{cc}
0 & 0 \\ 
0 & tv%
\end{array}%
\right) ,tw=\left( 
\begin{array}{cc}
0 & tsv \\ 
0 & 0%
\end{array}%
\right) , \\
stw &=&\left( 
\begin{array}{cc}
0 & 0 \\ 
0 & tsv%
\end{array}%
\right) ,tsw=\left( 
\begin{array}{cc}
0 & t^{2}v \\ 
0 & stv%
\end{array}%
\right) , \\
\left( sw\right) \left( tw\right) &=&0,t^{2}w=0,\left( tw\right) ^{2}=0\text{%
.}
\end{eqnarray*}%
Therefore, 
\begin{eqnarray*}
\left( sw\right) \left( tw\right) &=&0,\left( tw\right) \left( sw\right)
^{k}=\left( 
\begin{array}{cc}
0 & \left( tsv\right) \left( tv\right) ^{k} \\ 
0 & 0%
\end{array}%
\right) ,\left( tsw\right) \left( tw\right) =0, \\
\left( tw\right) \left( sw\right) ^{k} &=&\left( 
\begin{array}{cc}
0 & \left( tsv\right) \left( tv\right) ^{k} \\ 
0 & 0%
\end{array}%
\right) ,\left( tw\right) \left( tsw\right) ^{k}=\left( 
\begin{array}{cc}
0 & \left( tsv\right) \left( stv\right) ^{k} \\ 
0 & 0%
\end{array}%
\right) , \\
\left( tsw\right) \left( stw\right) ^{k} &=&\left( 
\begin{array}{cc}
0 & \left( t^{2}v\right) \left( tsv\right) ^{k} \\ 
0 & \left( stv\right) \left( tsv\right) ^{k}%
\end{array}%
\right) \text{, }\left( stw\right) \left( tsw\right) ^{k}=\left( 
\begin{array}{cc}
0 & 0 \\ 
0 & \left( tsv\right) \left( stv\right) ^{k}%
\end{array}%
\right) \text{.}
\end{eqnarray*}%
2. Let $w=tv(t,st)s=\left( 
\begin{array}{cc}
tv & 0 \\ 
sv & 0%
\end{array}%
\right) $. Then,%
\begin{eqnarray*}
sw &=&\left( 
\begin{array}{cc}
0 & 0 \\ 
tv & 0%
\end{array}%
\right) ,tw=\left( 
\begin{array}{cc}
tsv & 0 \\ 
0 & 0%
\end{array}%
\right) \\
stw &=&\left( 
\begin{array}{cc}
0 & 0 \\ 
tsv & 0%
\end{array}%
\right) ,tsw=\left( 
\begin{array}{cc}
t^{2}v & 0 \\ 
stv & 0%
\end{array}%
\right) \text{;}
\end{eqnarray*}%
\begin{eqnarray*}
t^{2}w &=&0,\left( tw\right) \left( sw\right) =0, \\
\left( sw\right) \left( tw\right) ^{k} &=&\left( 
\begin{array}{cc}
0 & 0 \\ 
\left( tv\right) \left( tsv\right) ^{k} & 0%
\end{array}%
\right) ,\left( tw\right) \left( tsw\right) ^{k}=\left( 
\begin{array}{cc}
\left( tsv\right) \left( t^{2}v\right) ^{k} & 0 \\ 
0 & 0%
\end{array}%
\right) , \\
stw\left( tsw\right) ^{k} &=&\left( 
\begin{array}{cc}
0 & 0 \\ 
\left( tsv\right) \left( t^{2}v\right) ^{k} & 0%
\end{array}%
\right) ,\left( tsw\right) \left( tw\right) ^{k}=\left( 
\begin{array}{cc}
t^{2}v\left( tsv\right) ^{k} & 0 \\ 
stv\left( tsv\right) ^{k} & 0%
\end{array}%
\right) , \\
\left( tsw\right) \left( stw\right) &=&0\text{.}
\end{eqnarray*}
\end{proof}

\begin{lemma}
The quotient of the ring $R$ by the ideal generated by $\left( tst\right)
^{3}$ has finite $\mathbb{Z}$-rank.
\end{lemma}

\begin{proof}
Let $u=tst$ and $N$ be the $2$-sided ideal in $R$ generated by $u^{3}$. We
need to consider only monomials which contain the subword $t^{2}$. Since $%
t^{3}=\left( st\right) ^{3}=0$, non-zero monomials which contain two
occurrences of $t^{2}$ have subwords of the form $t^{2}st^{2},t^{2}stst^{2}$%
. From the form $u^{3}=tst^{2}st^{2}st$, we conclude there are only a finite
number of monomials in $\frac{R}{N}$ with subwords $N+t^{2}st^{2}$. A long
enough non-zero monomial in $\frac{R}{N}$ will contain a subword $%
N+tstst^{2}stst^{2}stst^{2}stst=N+v^{4}$ where $v=tstst=\left( 
\begin{array}{cc}
0 & 0 \\ 
0 & st^{2}%
\end{array}%
\right) $; but as $v^{4}=0$, we are done.
\end{proof}

\begin{proposition}
(i) The set $S(w,4)=\left\{ 0\right\} $ for all $w\in \mu \left( R\right) $.%
\newline
(ii) The set $S(w,3)$ is contained in $L_{\infty }\left( R\right) $.\newline
(iii) $\mu \left( R\right) $ is nil of degree $5$.\newline
(iv) The semigroup $\mu \left( \frac{R}{L_{\infty }\left( R\right) }\right) $
is nil of degree $4$.
\end{proposition}

\begin{proof}
First we evaluate the set $S\left( w,k\right) $ for all monomials $w$ of
length at most $2$; therefore, we need to substitute only $w=1,t,ts$. The
set of values produced is 
\begin{equation*}
M\left( k\right) =\left\{ 
\begin{array}{c}
ts^{k},st^{k},t\left( ts\right) ^{k},\left( ts\right) \left( t^{2}\right)
^{k},t^{2}\left( st\right) ^{k}, \\ 
\left( tst\right) \left( st^{2}\right) ^{k},\left( sts\right) \left(
t^{2}s\right) ^{k},\left( t^{2}s\right) \left( ts\right) ^{2k}%
\end{array}%
\right\} \text{.}
\end{equation*}%
More specifically, we find 
\begin{eqnarray*}
M\left( 1\right) &=&\left\{
ts,st,t^{2}s,tst^{2},t^{2}st,tstst^{2},stst^{2}s,t^{2}ststs\right\} \text{,}
\\
M\left( 2\right) &=&\left\{ st^{2},t\left( ts\right) ^{2},t^{2}\left(
st\right) ^{2},\left( tst\right) \left( st^{2}\right) ^{2},\left( sts\right)
\left( t^{2}s\right) ^{2}\right\} \text{,} \\
M\left( 3\right) &=&\left\{ t\left( ts\right) ^{3},\left( tst\right) \left(
st^{2}\right) ^{3},\left( sts\right) \left( t^{2}s\right) ^{3}\right\} \text{%
,} \\
M\left( k\right) &=&\left\{ 0\right\} \text{ for }k\geq 4\text{.}
\end{eqnarray*}

(i) To prove $S\left( w,4\right) =\left\{ 0\right\} $, we proceed by
induction on the length of $w$. The fact $M\left( 4\right) =\left\{
0\right\} $ takes care of the first step of the induction. Suppose $w$ has
length greater than $2$, then by the previous proposition, $S\left(
w,4\right) $ is contained in $\cup \left\{ S(v,k)\mid \left\vert
v\right\vert <\left\vert w\right\vert \right\} $ and thus we are done by
induction.

(ii) We note that $\left( ts\right) ^{3},\left( st^{2}\right) ^{3},\left(
t^{2}s\right) ^{3}\in L_{2}\left( R\right) $ and therefore $M\left( 3\right)
\ $is contained in $L_{2}\left( R\right) $. By induction on the length of $w$
and by using 
\begin{equation*}
\left( 
\begin{array}{cc}
L_{i}\left( R\right) & L_{i}\left( R\right) \\ 
L_{i}\left( R\right) & L_{i}\left( R\right)%
\end{array}%
\right) \cap R\leq L_{i+1}\left( R\right)
\end{equation*}%
we reach $S\left( w,3\right) \leq L_{\infty }\left( R\right) $.

(iii) To prove that $w\in \mu \left( R\right) $ is nil of degree at most $5$%
, we proceed by induction on $\left\vert w\right\vert $ and use part (i)
above. To prove that the degree is exactly $5$, it is sufficient to check
that $\delta \left( tst^{2}s\right) =5$.

(iv) To prove that $\overline{w}\in $ $\overline{\mu \left( R\right) }$ is
nil of degree at most $4$, we proceed by induction on $\left\vert
w\right\vert $ and use part (ii) above. Suppose the nil degree is less than $%
4$ then $u=\left( tst\right) ^{3}\in L_{\infty }\left( R\right) $ and by
Lemma 9, $\overline{R}$ has finite $\mathbb{Z}$-rank, which contradicts
Proposition 2.
\end{proof}

We note here that $\mu \left( R\right) $ is not a maximal nil multiplicative
semigroup. Indeed, for $u=s+t^{2}$, we easily check that $\delta \left(
u\right) =8$ and 
\begin{equation*}
u.\mu \left( R\right) ,\mu \left( R\right) .u\leq \mu \left( R\right) \text{;%
}
\end{equation*}
thus, the multiplicative semigroup $\left\langle u,s,t\right\rangle $ is nil
of degree $8$.

\subsection{Presentation of $R$}

Recall the elements $u_{1}=\sigma ^{2},$ $u_{2i}=\tau .u_{2i-1}^{\lambda },$ 
$u_{2i+1}=u_{2i}^{\lambda }$ $\left( i\geq 1\right) $ in the free ring $%
F=\left\langle \sigma ,\tau \right\rangle $ and the endomorphism of $F$
defined by $\lambda :\sigma \rightarrow \tau ,$ $\tau \rightarrow \sigma
\tau $.

\begin{proposition}
The ring $R$ has the presentation 
\begin{equation*}
\left\langle \sigma ,\tau \mid u_{i}=0\text{ }\left( i\geq 1\right)
\right\rangle \text{.}
\end{equation*}
\end{proposition}

\begin{proof}
(1) By Proposition 1, if $w\left( s,t\right) =0$ is a minimal relation then $%
w\left( s,t\right) $ is a monomial.

(2) Let $w\left( s,t\right) =0$ be minimal and $w\left( s,t\right) =v\left(
s,t\right) s$. Then, by Lemma 6, $v\left( s,t\right) \in lan\left( s\right) =%
\widehat{R}s$ and thus, $w\left( s,t\right) =s^{2}$.

(3) Let $w\left( s,t\right) =0$ and $w\left( s,t\right) =v\left( s,t\right)
t $. Then, by Lemma 3, $w\left( s,t\right) =stw^{\prime }\left( t,st\right) $
or $tw^{\prime }\left( t,st\right) $.

(3.1) If $w\left( s,t\right) =stw^{\prime }\left( t,st\right) $ then $%
w\left( s,t\right) ^{\varphi }=tw^{\prime }\left( s,t\right) =0$.

On the other hand, $tw^{\prime }\left( s,t\right) =0$ $\Rightarrow $ $%
stw^{\prime }\left( t,st\right) =\left( tw^{\prime }\left( s,t\right)
\right) ^{\lambda }=0$.

We assert: $tw^{\prime }\left( s,t\right) =0$ not minimal $\Rightarrow $ $%
w\left( s,t\right) =stw^{\prime }\left( t,st\right) =0$ not minimal.

For suppose, $tw^{\prime }\left( s,t\right) =tu_{1}u_{2}$ and $u_{2}=0$. If $%
u_{2}=tx_{2}$ then $w\left( s,t\right) =st\left( u_{1}\right) ^{\lambda
}\left( u_{2}\right) ^{\lambda }$ and $\left( u_{2}\right) ^{\lambda
}=st\left( x_{2}\right) ^{\lambda }=0$ and not minimal.

Similarly, if $u_{1}=r_{1}tr_{2}$.

As $u_{2}=sx_{2}$, we are left with the possibility $u_{1}$ empty; that is, $%
w^{\prime }\left( s,t\right) =u_{2}$. But then, $tw^{\prime }\left(
t,st\right) =0$ and $w\left( s,t\right) =stw^{\prime }\left( t,st\right) $
is not minimal.

(3.2) If $w\left( s,t\right) =tw^{\prime }\left( t,st\right) $ then $%
w^{\prime }\left( s,t\right) =0$.

We assert: $w^{\prime }\left( s,t\right) =0$ not minimal $\Rightarrow $ $%
w\left( s,t\right) =tw^{\prime }\left( t,st\right) =0$ not minimal.

Suppose $w^{\prime }\left( s,t\right) =u_{1}u_{2}u_{3}$ where $u_{2}=0$
minimal and $u_{1}$ or $u_{3}$ non-empty. Furthermore suppose $w\left(
s,t\right) =0$ minimal.

We have, $w\left( s,t\right) =t\left( u_{1}\right) ^{\lambda }\left(
u_{2}\right) ^{\lambda }\left( u_{3}\right) ^{\lambda }$.

If $u_{1}$ is empty then $w\left( s,t\right) =t\left( u_{2}\right) ^{\lambda
}\left( u_{3}\right) ^{\lambda }$ and $t\left( u_{2}\right) ^{\lambda
}=0,u_{3}$ empty; contradiction.

We may assume $u_{1}\not=0$.

If $u_{2}=tu_{2}^{\prime }$ then $\left( u_{2}\right) ^{\lambda }=0$; thus, $%
u_{2}=su_{2}^{\prime }$.

If $u_{1}=u_{1}^{\prime }tu_{1}^{\prime \prime }$ then $w\left( s,t\right)
=t\left( u_{1}^{\prime }\right) ^{\lambda }\left( tu_{1}^{\prime \prime
}u_{2}\right) ^{\lambda }\left( u_{3}\right) ^{\lambda }$ and $\left(
tu_{1}^{\prime \prime }u_{2}\right) ^{\lambda }=0$; contradiction.

Then, $u_{1}=s$ and so, $w^{\prime }\left( s,t\right) =s\left(
su_{2}^{\prime }\right) u_{3}=$.$s^{2}u_{4}$. As $s^{2}=0$, we are done.

Hence, the set of minimal monomial relators is 
\begin{equation*}
\left\{ 
\begin{array}{c}
u_{1}\left( s,t\right) =s^{2},\text{ }u_{2i}\left( s,t\right)
=t.u_{2i-1}\left( s,t\right) ^{\lambda }, \\ 
u_{2i+1}\left( s,t\right) =u_{2i}\left( s,t\right) ^{\lambda }\text{ }\left(
i\geq 1\right)%
\end{array}%
\right\}
\end{equation*}%
where $\lambda :s\rightarrow t,t\rightarrow st$.
\end{proof}

\subsection{Growth of $R$}

Let $F\left( n\right) $ be the number of non-zero elements of $\mu \left(
R\right) $ of length at most $n$. This function is the growth function for $%
\mu \left( R\right) $ and also of $R$. Let $V\left( n\right) $ to be the set
of non-zero elements of $\mu \left( R\right) $ of length $n$ and let $%
f\left( n\right) =\mid V\left( n\right) \mid $. Computer enumeration
provides the following values%
\begin{eqnarray*}
f(n) &=&2,3,4,5,7,8,9,11,13,15, \\
&&16,17,19,21,24,27,29,31,32,33, \\
&&35,37,40,43,46,50,53,56,59,61, \\
&&63,64,65,67,69,72,75,78,82,85, \\
&&89,94,98,103,107,110,114,117,120,123, \\
&&125,127,128,129,131,133,136,139,142,146,...\text{.}
\end{eqnarray*}

We note the occurrence of powers of $2^{i}$ $\left( 1\leq i\leq 7\right) $
in this list. This suggests that $f\left( a_{i}-2\right) =2^{i-3}$ for $%
i\geq 4$ where $a_{i}$ is the $i$th term of the Fibonacci sequence ($%
a_{0}=0,a_{1}=1,...,a_{i}=a_{i-1}+a_{i-2}$ for $i\geq 2$). The proof of this
observation is a lengthy argument. The appearance of the Fibonacci sequence
can be explained in part by the fact that the linearized form of the
endomorphism $\lambda :\sigma \rightarrow \tau ,\tau \rightarrow \sigma \tau 
$ which produces the relations of $R$, has the matrix representation%
\begin{equation*}
\lambda ^{\ast }=\left( 
\begin{array}{cc}
0 & 1 \\ 
1 & 1%
\end{array}%
\right) ,\left( \lambda ^{\ast }\right) ^{i}=\left( 
\begin{array}{cc}
a_{i-1} & a_{i} \\ 
a_{i} & a_{i+1}%
\end{array}%
\right)
\end{equation*}%
for all $i$.

To start we show

\begin{lemma}
Suppose $f$ is strictly monotone increasing and that $f\left( a_{i}-2\right)
=2^{i-3}$ for all $i\geq 4$. Then, 
\begin{equation*}
\frac{1}{8}n^{c}<f\left( n\right) <2n^{c}
\end{equation*}%
where $c=\frac{\log \left( 2\right) }{\log \left( \alpha \right) },\alpha =%
\frac{1+\sqrt{5}}{2}$.
\end{lemma}

\begin{proof}
First, we recall that , 
\begin{equation*}
\alpha ^{i-2}<a_{i}<\alpha ^{i-1}
\end{equation*}%
holds for all $i\geq 3$, where $\alpha =\frac{1+\sqrt{5}}{2}$.

Let $n$ be a length of a nonzero monomial and let $k$ be such that $%
a_{k}\leq n<a_{k+1}$. Then,%
\begin{eqnarray*}
\alpha ^{k-2} &<&n<\alpha ^{k}, \\
\log _{\alpha }\left( n\right) &<&k<\log _{\alpha }\left( n\right) +2
\end{eqnarray*}%
Since $f$ is a strictly monotone increasing function, we have%
\begin{equation*}
f\left( a_{k}\right) \leq f\left( n\right) <f\left( a_{k+1}\right) \text{,}
\end{equation*}%
\begin{equation*}
2^{\log _{\alpha }\left( n\right) -3}<2^{k-3}<f\left( n\right)
<2^{k-1}<2^{\log _{\alpha }\left( n\right) +1}
\end{equation*}%
and thus%
\begin{equation*}
\frac{1}{8}n^{c}<f\left( n\right) <2n^{c}\text{.}
\end{equation*}
\end{proof}

We will reduce the problem of the growth of $V\left( n\right) $ to the
growth of%
\begin{equation*}
W\left( n\right) =\left\{ v\in V\left( n\right) \mid
sv\not=0,tv\not=0\right\} \text{,}
\end{equation*}%
by showing that 
\begin{equation*}
g\left( n\right) =\mid W\left( n\right) \mid =f\left( n+1\right) -f\left(
n\right) \text{,}
\end{equation*}%
the discrete derivative of $f$.

Toward this end we let $U$ be the set of monomial relations of $R$,

\begin{eqnarray*}
s\backslash U &=&\left\{ v\mid sv\in U\right\} ,t\backslash U=\left\{ v\mid
tv\in U\right\} , \\
_{s}V\left( n\right) &=&\left\{ v\in V\left( n\right) \mid sv=0\right\}
,_{t}V\left( n\right) =\left\{ v\in V\left( n\right) \mid tv=0\right\} , \\
_{\overline{s}}V\left( n\right) &=&\left\{ v\in V\left( n\right) \mid
sv\not=0\right\} ,_{\overline{t}}V\left( n\right) =\left\{ v\in V\left(
n\right) \mid tv\not=0\right\} \text{.}
\end{eqnarray*}

\begin{lemma}
$\mid W\left( n\right) \mid =g\left( n\right) $ for all $n\geq 1$.
\end{lemma}

\begin{proof}
As $ran\left( R\right) =0$, we have 
\begin{equation*}
_{s}V\left( n\right) \leq \text{ }_{\overline{t}}V\left( n\right) ,\text{ }%
_{t}V\left( n\right) \leq \text{ }_{\overline{s}}V\left( n\right) .
\end{equation*}%
Then, $V\left( n\right) $ is partitioned as%
\begin{eqnarray*}
V\left( n\right) &=&\text{ }_{s}V\left( n\right) \cup \text{ }_{t}V\left(
n\right) \cup W\left( n\right) \text{,} \\
sV\left( n\right) &=&s\left( _{t}V\left( n\right) \right) \cup sW\left(
n\right) , \\
tV\left( n\right) &=&t\left( _{s}V\left( n\right) \right) \cup tW\left(
n\right) \text{.}
\end{eqnarray*}%
\begin{eqnarray*}
V\left( n\right) &=&_{s}V\left( n\right) \cup \text{ }_{\overline{s}}V\left(
n\right) =\text{ }_{t}V\left( n\right) \cup \text{ }_{\overline{t}}V\left(
n\right) , \\
&=&\text{ }_{\overline{s}}V\left( n\right) \cup \text{ }_{\overline{t}%
}V\left( n\right) \text{.}
\end{eqnarray*}%
We have 
\begin{eqnarray*}
V\left( n+1\right) &=&sV\left( n\right) \cup tV\left( n\right) =s\left(
_{t}V\left( n\right) \right) \cup sW\left( n\right) \cup t\left( _{s}V\left(
n\right) \right) \cup tW\left( n\right) , \\
&=&s\left( _{\overline{s}}V\left( n\right) \right) \cup s\left( _{\overline{t%
}}V\left( n\right) \right) \cup t\left( _{\overline{s}}V\left( n\right)
\right) \cup t\left( _{\overline{t}}V\left( n\right) \right) \text{.}
\end{eqnarray*}%
From the partition%
\begin{equation*}
V\left( n+1\right) =sV\left( n\right) \cup tV\left( n\right) =\left( s\left(
_{t}V\left( n\right) \right) \cup sW\left( n\right) \right) \cup \left(
t\left( _{s}V\left( n\right) \right) \cup tW\left( n\right) \right) \text{,}
\end{equation*}%
we conclude%
\begin{eqnarray*}
&\mid &V\left( n+1\right) \mid =\mid _{t}V\left( n\right) \mid +\mid W\left(
n\right) \mid +\mid _{s}V\left( n\right) \mid +\mid W\left( n\right) \mid \\
&=&\left( \mid _{t}V\left( n\right) \mid +\mid _{s}V\left( n\right) \mid
+\mid W\left( n\right) \mid \right) +\mid W\left( n\right) \mid \\
&=&\mid V\left( n\right) \mid +\mid W\left( n\right) \mid \text{.}
\end{eqnarray*}
\end{proof}

Next we show that there are two operations acting on $W$.

\begin{proposition}
(i)Let $W_{s},W_{t}$ be the subsets of monomials in $W$ which end in $s,t$,
respectively. Then\newline
$W_{s}=W_{t}s,W_{t}=\left( sW\right) ^{\lambda }$ \newline
(ii) The function $f$ is strictly monotone increasing.\newline
\end{proposition}

\begin{proof}
(i) To start with, $t\in W$ since $s.t\not=0\not=t.t$.

Let $w=w^{\prime }t\in W\left( n\right) $ then $ws=w^{\prime }ts\not=0$ and $%
sws=sw^{\prime }ts,tws=tw^{\prime }ts$ are both non-zero and so, $ws\in
W\left( n+1\right) $

Let $w=w^{\prime }s\in W(n)$. Then $wt\in W(n+1)$ iff $wt\not=s\backslash
u,t\backslash u$.

It $w=tsqt\in W\left( n+1\right) $ then $sqt=v^{\lambda }$, $w=tv^{\lambda }$
and $v=tv^{\prime }$; for if $v=sv^{\prime }$ then $v^{\lambda }=t\ast $ and 
$w=t^{2}\ast $, a contradiction.

Suppose $sv=0$; then $v=\left( s\backslash u_{2n}\right) r=u^{\prime }r$.
However,%
\begin{eqnarray*}
u_{2n} &=&su^{\prime },u_{2n+1}=t\left( u_{2n}\right) ^{\lambda
}=t^{2}\left( u^{\prime }\right) ^{\lambda }=0, \\
w &=&tv^{\lambda }=t\left( u^{\prime }r\right) ^{\lambda }=t\left( u^{\prime
}\right) ^{\lambda }r^{\lambda }andtw=0\text{.}
\end{eqnarray*}

Suppose $tv=0$; then $v=\left( t\backslash u_{2n-1}\right) r=u^{\prime }r$.
However,%
\begin{eqnarray*}
u_{2n-1} &=&tu^{\prime },u_{2n}=st\left( u^{\prime }\right) ^{\lambda }=0, \\
w &=&tv^{\lambda }=t\left( u^{\prime }\right) ^{\lambda }r^{\lambda }
\end{eqnarray*}
and $sw=0$.

On the other hand suppose $v=v\left( s,t\right) \in W\left( k\right) $ of
type $\left( k_{s},k_{t}\right) $. Then, $\left( sv\left( s,t\right) \right)
^{\lambda }=t.v\left( t,st\right) =\left( 
\begin{array}{cc}
0 & t.v\left( s,t\right) \\ 
0 & s.v\left( s,t\right)%
\end{array}%
\right) $ and

\begin{equation*}
st.v\left( t,st\right) =\left( 
\begin{array}{cc}
0 & 0 \\ 
0 & t.v\left( s,t\right)%
\end{array}%
\right) ,t^{2}.v\left( t,st\right) =\left( 
\begin{array}{cc}
0 & ts.v(s,t) \\ 
0 & 0%
\end{array}%
\right)
\end{equation*}%
both non-zero and $\left( sv\left( s,t\right) \right) ^{\lambda }$ is of
type $\left( k_{t},1+k_{s}+k_{t}\right) =\left( k_{t},1+k\right) $ has
length $1+n=1+k_{s}+2k_{t}$.

(ii) Since $ran\left( R\right) =0$, we conclude that $f$ is monotone
increasing.
\end{proof}

\subsubsection{Word types in $W$}

We have seen that the following two transformations act on $W\cup \left\{
0\right\} $:%
\begin{equation*}
\sigma :w\rightarrow ws,\kappa :w\rightarrow \left( sw\right) ^{\lambda }%
\text{,}
\end{equation*}%
where $\sigma ^{2}=0$ and where $\left\{ \left( w\right) \kappa ^{n}\mid
n\geq 0\right\} $ is an infinite set for all $w\in W$. Given positive
integers $m,n_{1},n_{2}...,n_{m}$, let $q=q\left(
n_{1},n_{2}...,n_{m}\right) =\kappa ^{n_{1}}\sigma \kappa ^{n_{2}}...\sigma
\kappa ^{n_{m}}$.

Let $\zeta $ be the set of transformations consisting of $e$ (the trivial
operator)$,\sigma $ and $q,\sigma q,q\sigma ,\sigma q\sigma $ for all
possible $q$. Then, $W=\left( t\right) \zeta $.

Given $w\in W$, let $c$ denote the $s$-length, $d$ the $t$-length of $w$ and
call $\left\Vert w\right\Vert =\left( c,d\right) $ the type of $w$. Then $%
\left( W\right) \subset \mathbb{N\times N}$ and 
\begin{equation*}
\left\Vert \left( w\right) \sigma \right\Vert =\left( c+1,d\right)
,\left\Vert \left( w\right) \kappa \right\Vert =\left( d,c+d+1\right) \text{.%
}
\end{equation*}%
Corresponding to $\sigma ,\lambda ,\kappa $, we define operations $\sigma
^{\ast },\lambda ^{\ast },\kappa ^{\ast }$ on $\mathbb{Z\times Z}$:%
\begin{eqnarray*}
\sigma ^{\ast } &:&\left( c,d\right) \rightarrow \left( c+1,d\right) \text{
translation;} \\
\lambda ^{\ast } &:&\left( c,d\right) \rightarrow \left( d,c+d\right) \text{,%
} \\
\kappa ^{\ast }\left( =\sigma ^{\ast }\lambda ^{\ast }\right) &:&\left(
c,d\right) \rightarrow \left( d,c+d+1\right) \text{.}
\end{eqnarray*}%
Then, $\sigma ^{\ast },\lambda ^{\ast }$ are invertible transformations. Let 
$\rho ^{\ast }$ be the group generated by $\sigma ^{\ast },\lambda ^{\ast }$%
. We verify 
\begin{equation*}
\sigma ^{\ast }\left( \sigma ^{\ast }\right) ^{\left( \lambda ^{\ast
}\right) ^{-1}}=\left( \sigma ^{\ast }\right) ^{\lambda ^{\ast }}\text{,}
\end{equation*}%
\begin{equation*}
\left\langle \sigma ^{\ast },\left( \sigma ^{\ast }\right) ^{\lambda ^{\ast
}}\right\rangle \cong \mathbb{Z\times Z}
\end{equation*}%
and $\rho ^{\ast }$ is metabelian group, an extension of the normal closure
of $\left\langle \sigma ^{\ast }\right\rangle $ by $\left\langle \lambda
^{\ast }\right\rangle $.

We write $\lambda ^{\ast }=\left( 
\begin{array}{cc}
0 & 1 \\ 
1 & 1%
\end{array}%
\right) $ with respect to the canonical basis of $\mathbb{Z\times Z}$, with
action on the right hand side. Then as noted before, for any integer $i$ we
have 
\begin{equation*}
\left( \lambda ^{\ast }\right) ^{i}=\left( 
\begin{array}{cc}
a_{i-1} & a_{i} \\ 
a_{i} & a_{i+1}%
\end{array}%
\right)
\end{equation*}%
where $a_{i}$ is the $i$th term of the extended Fibonacci sequence defined
by $a_{0}=0,a_{1}=1,...,a_{i}=a_{i-1}+a_{i-2}$ for $i\geq 2$ and $%
a_{-i}=\left( -1\right) ^{i+1}a_{i}$ for $i\geq 1$. We gather together in
the following lemma some facts about the Fibonacci numbers which we will
need in the sequel.

\begin{lemma}
(Some Fibonacci facts) The following equations hold 
\begin{equation*}
a_{1}+a_{2}+...+a_{n}=a_{n+2}-1,
\end{equation*}%
\begin{equation*}
a_{1}+a_{3}+a_{5}+...+a_{2n+1}=a_{2n+2}
\end{equation*}%
for all integers $n\geq 1$;%
\begin{equation*}
a_{m}a_{n}+a_{m-1}a_{n-1}=a_{m+n-1}
\end{equation*}%
for all integers $m,n$.
\end{lemma}

We consider within $\rho ^{\ast }$ the set of words $\zeta ^{\ast }$
corresponding to $\zeta $.

For positive integers $m,n_{1},n_{2}...,n_{m}$ and $n=n_{1}+n_{2}+...+n_{m}$%
,there are $4$ word forms in $\zeta ^{\ast }$: 
\begin{eqnarray*}
\kappa ^{\ast n_{1}}\sigma ^{\ast }...\kappa ^{\ast n_{m}}\sigma ^{\ast }
&=&\left( \sigma ^{\ast }\right) ^{\left( \lambda ^{\ast }\right)
^{-n_{1}}}\left( \sigma ^{\ast }\right) ^{\left( \lambda ^{\ast }\right)
^{-\left( n_{1}+n_{2}\right) }}...\left( \sigma ^{\ast }\right) ^{\left(
\lambda ^{\ast }\right) ^{-n}}\kappa ^{\ast n}, \\
\sigma ^{\ast }\kappa ^{\ast n_{1}}\sigma ^{\ast }...\kappa ^{\ast
n_{m}}\sigma ^{\ast } &=&\sigma ^{\ast }\left( \sigma ^{\ast }\right)
^{\left( \lambda ^{\ast }\right) ^{-n_{1}}}\left( \sigma ^{\ast }\right)
^{\left( \lambda ^{\ast }\right) ^{-\left( n_{1}+n_{2}\right) }}...\left(
\sigma ^{\ast }\right) ^{\left( \lambda ^{\ast }\right) ^{-n}}\kappa ^{\ast
n}, \\
\kappa ^{\ast n_{1}}\sigma ^{\ast }...\sigma ^{\ast }\kappa ^{\ast n_{m}}
&=&\left( \sigma ^{\ast }\right) ^{\left( \lambda ^{\ast }\right)
^{-n_{1}}}\left( \sigma ^{\ast }\right) ^{\left( \lambda ^{\ast }\right)
^{-\left( n_{1}+n_{2}\right) }}...\left( \sigma ^{\ast }\right) ^{\left(
\lambda ^{\ast }\right) ^{-\left( n_{1}+n_{2}+...+n_{m-1}\right) }}\kappa
^{\ast n}, \\
\sigma ^{\ast }\kappa ^{\ast n_{1}}\sigma ^{\ast }...\sigma ^{\ast }\kappa
^{\ast n_{m}} &=&\sigma ^{\ast }\left( \sigma ^{\ast }\right) ^{\left(
\lambda ^{\ast }\right) ^{-n_{1}}}\left( \sigma ^{\ast }\right) ^{\left(
\lambda ^{\ast }\right) ^{-\left( n_{1}+n_{2}\right) }}...\left( \sigma
^{\ast }\right) ^{\left( \lambda ^{\ast }\right) ^{-\left(
n_{1}+n_{2}+...+n_{m-1}\right) }}\kappa ^{\ast n}\text{.}
\end{eqnarray*}%
We need to introduce numerical functions $\Delta _{1}\left(
n:n_{1},...,n_{m}\right) ,\Delta _{2}\left( n:n_{1},...,n_{m}\right) $
defined by:%
\begin{equation*}
\Delta _{1}\left( n:n_{1},...,n_{m}\right) =\Delta _{2}\left(
n:n_{1},...,n_{m}\right) =0\text{ for }m=1
\end{equation*}%
and 
\begin{eqnarray*}
\Delta _{1}\left( n:n_{1},...,n_{m}\right) &=&a_{n-\left( n_{1}+1\right)
}+a_{n-\left( n_{1}+n_{2}+1\right) }+...+a_{n-\left(
n_{1}+n_{2}+...+n_{m-1}+1\right) }, \\
\Delta _{2}\left( n:n_{1},...,n_{m}\right) &=&a_{n-n_{1}}+a_{n-\left(
n_{1}+n_{2}\right) }+...+a_{n-\left( n_{1}+n_{2}+...+n_{m-1}\right) }\text{.}
\end{eqnarray*}%
for $m>1$.

\begin{lemma}
The functions $\Delta _{1},\Delta _{2}$ satisfy the following properties:$%
\newline
$(i) $\Delta _{1}\left( n:n_{1},...,n_{m}\right) \leq \Delta _{1}\left(
n:1,...,1\right) =\sum_{n-2}a_{i}=a_{n}-1$;$\newline
$(ii) $\Delta _{2}\left( n:n_{1},...,n_{m}\right) \leq \Delta _{2}\left(
n:1,...,1\right) =\sum_{n}a_{i}=a_{n+1}-1$;\newline
(iii) $\Delta _{1}\left( n:n_{1},...,n_{m}\right) =0\Leftrightarrow m=1$ or $%
m=2,n_{1}=n-1,n_{2}=1$;\newline
(iv) $\Delta _{2}\left( n:n_{1},...,n_{m}\right) =0\Leftrightarrow m=1$ or $%
m=2,n=n_{1}$;\newline
(v) $\Delta _{2}\left( n:n_{1},...,n_{m}\right) =1\Leftrightarrow
m=2,n_{1}=n-1,n_{2}=1$;\newline
(vi) for $j=1,2,\newline
\Delta _{j}\left( n:n_{1},...,n_{m}\right) =\Delta _{j}\left(
n:1,...,1\right) \Leftrightarrow n_{i}=1$ for all $i$;\newline
(vii)%
\begin{eqnarray*}
&&\Delta \left( n:n_{1},...,n_{m}\right) =\Delta _{1}\left(
n:n_{1},...,n_{m}\right) +\Delta _{2}\left( n:n_{1},...,n_{m}\right) \\
&=&a_{n-n_{1}+1}+a_{n-\left( n_{1}+n_{2}\right) +1}+...+a_{n-\left(
n_{1}+n_{2}+...+n_{m-1}\right) +1}\leq a_{n+2}-2\text{;}
\end{eqnarray*}%
(viii)%
\begin{equation*}
\Delta _{2}-\Delta _{1}=a_{n-\left( n_{1}+2\right) }+a_{n-\left(
n_{1}+n_{2}+2\right) }+...+a_{n-\left( n_{1}+n_{2}+...+n_{m-1}+\delta
\right) }
\end{equation*}%
where $\delta =1$ if $n_{m}=1$ and $\delta =2$ otherwise;\newline
(ix) $\Delta _{2}-\Delta _{1}\leq a_{n-n_{1}}-1$;$\newline
$(x) $\Delta _{2}=\Delta _{1}$ iff $\Delta _{2}=\Delta _{1}=0$ or $%
m=2,n_{1}=n-2,n_{2}=2$.
\end{lemma}

In the next two lemmas we determine the image of a word type $\left(
c,d\right) $ under the application of elements of the group $\zeta ^{\ast }$.

\begin{lemma}
Let $j\geq 0$. Then,%
\begin{eqnarray*}
\left( c,d\right) .\left( \sigma ^{\ast }\right) ^{\left( \lambda ^{\ast
}\right) ^{-j}} &=&(c+a_{-j-1},d+a_{-j}), \\
\left( c,d\right) .\left( \kappa ^{\ast }\right) ^{j} &=&\left(
a_{j-1}c+a_{j}d+a_{j+1}-1,a_{j}c+a_{j+1}d+a_{j+2}-1\right) , \\
\left( 0,1\right) \left( \kappa ^{\ast }\right) ^{j} &=&\left(
a_{j+2}-1,a_{j+3}-1\right) \text{.}
\end{eqnarray*}
\end{lemma}

\begin{proof}
(i) 
\begin{eqnarray*}
\left( c,d\right) \left( \lambda ^{\ast }\right) ^{j}\sigma ^{\ast }\left(
\lambda ^{\ast }\right) ^{-j} &=&\left( c,d\right) \left( 
\begin{array}{cc}
a_{j-1} & a_{j} \\ 
a_{j} & a_{j+1}%
\end{array}%
\right) \sigma ^{\ast }\left( \lambda ^{\ast }\right) ^{-j} \\
&=&\left( ca_{j-1}+da_{j}+1,ca_{j}+da_{j+1}\right) \left( 
\begin{array}{cc}
a_{-j-1} & a_{-j} \\ 
a_{-j} & a_{-j+1}%
\end{array}%
\right) \\
&=&(\left( ca_{j-1}+da_{j}+1\right) a_{-j-1}+\left( ca_{j}+da_{j+1}\right)
a_{-j}, \\
&&\left( ca_{j-1}+da_{j}+1\right) a_{-j}+\left( ca_{j}+da_{j+1}\right)
a_{-j+1}) \\
&=&(ca_{j-1}a_{-j-1}+da_{j}a_{-j-1}+a_{-j-1}+ca_{j}a_{-j}+da_{j+1}a_{-j}, \\
&&ca_{j-1}a_{-j}+da_{j}a_{-j}+a_{-j}+ca_{j}a_{-j+1}+da_{j+1}a_{-j+1}) \\
&=&(c\left( a_{j-1}a_{-j-1}+a_{j}a_{-j}\right) +d\left(
a_{j}a_{-j-1}+a_{j+1}a_{-j}\right) +a_{-j-1}, \\
&&c\left( a_{j-1}a_{-j}+a_{j}a_{-j+1}\right) +d\left(
a_{j}a_{-j}+a_{j+1}a_{-j+1}\right) +a_{-j}) \\
&=&(c+a_{-j-1},d+a_{-j})\text{.}
\end{eqnarray*}

(ii) Since $\left( \kappa ^{\ast }\right) ^{j}=\sigma ^{\ast }\left( \sigma
^{\ast }\right) ^{\left( \lambda ^{\ast }\right) ^{-1}}...\left( \sigma
^{\ast }\right) ^{\left( \lambda ^{\ast }\right) ^{-(j-1)}}\left( \lambda
^{\ast }\right) ^{j}$, we find%
\begin{equation*}
\left( c,d\right) \left( \kappa ^{\ast }\right)
^{j}=(c+s_{-j},d+s_{-j+1})\left( \lambda ^{\ast }\right) ^{j}
\end{equation*}

where%
\begin{eqnarray*}
s_{-j} &=&a_{-1}+a_{-2}+...+a_{-j}, \\
s_{-j+1} &=&a_{0}+a_{-1}+...+a_{-j+1}\text{.}
\end{eqnarray*}

Thus,%
\begin{eqnarray*}
\left( c,d\right) \left( \kappa ^{\ast }\right) ^{j}
&=&(c+s_{-j},d+s_{-j+1})\left( 
\begin{array}{cc}
a_{j-1} & a_{j} \\ 
a_{j} & a_{j+1}%
\end{array}%
\right) \\
&=&(ca_{j-1}+da_{j}+\left( s_{-j}a_{j-1}+s_{-j+1}a_{j}\right)
,ca_{j}+da_{j+1}+\left( s_{-j}a_{j}+s_{-j+1}a_{j+1}\right) ) \\
&=&(ca_{j-1}+da_{j}+a_{j+1}-1,ca_{j}+da_{j+1}+a_{j+2}-1)\text{.}
\end{eqnarray*}

(iii) In particular,%
\begin{eqnarray*}
\left( 0,1\right) \left( \kappa ^{\ast }\right) ^{j} &=&\left(
a_{j}+a_{j+1}-1,a_{j+1}+a_{j+2}-1\right) \\
&=&\left( a_{j+2}-1,a_{j+3}-1\right) \text{.}
\end{eqnarray*}
\end{proof}

Next, we derive the following more general formulae.

\begin{lemma}
\begin{eqnarray*}
&&\text{(i) }\left( c,d\right) \left( \kappa ^{\ast }\right) ^{n_{1}}\sigma
^{\ast }...\left( \kappa ^{\ast }\right) ^{n_{m}}\sigma ^{\ast } \\
&=&(a_{n-1}c+a_{n}d+a_{n+1}+\Delta _{1}\left( n:n_{1},...,n_{m}\right) , \\
&&a_{n}c+a_{n+1}d+a_{n+2}+\Delta _{2}\left( n:n_{1},...,n_{m}\right) -1)%
\text{;}
\end{eqnarray*}%
\begin{eqnarray*}
&&\text{(ii) }\left( c,d\right) \left( \kappa ^{\ast }\right) ^{n_{1}}\sigma
^{\ast }...\left( \kappa ^{\ast }\right) ^{n_{m}} \\
&=&(a_{n-1}c+a_{n}d+a_{n+1}+\Delta _{1}\left( n:n_{1},...,n_{m}\right) -1, \\
&&a_{n}c+a_{n+1}d+a_{n+2}+\Delta _{2}\left( n:n_{1},...,n_{m}\right) -1)%
\text{;}
\end{eqnarray*}%
\begin{eqnarray*}
&&\text{(iii) }\left( c,d\right) \sigma ^{\ast }\left( \kappa ^{\ast
}\right) ^{n_{1}}\sigma ^{\ast }...\left( \kappa ^{\ast }\right)
^{n_{m}}\sigma ^{\ast } \\
&=&(a_{n-1}\left( c+1\right) +a_{n}d+a_{n+1}+\Delta _{1}\left(
n:n_{1},...,n_{m}\right) , \\
&&a_{n}\left( c+1\right) +a_{n+1}d+a_{n+2}+\Delta _{2}\left(
n:n_{1},...,n_{m}\right) -1)\text{;}
\end{eqnarray*}%
\begin{eqnarray*}
&&\text{(iv) }\left( c,d\right) \sigma ^{\ast }\left( \kappa ^{\ast }\right)
^{n_{1}}\sigma ^{\ast }...\left( \kappa ^{\ast }\right) ^{n_{m}} \\
&=&(a_{n-1}\left( c+1\right) +a_{n}d+a_{n+1}+\Delta _{1}\left(
n:n_{1},...,n_{m}\right) -1, \\
&&a_{n}\left( c+1\right) +a_{n+1}d+a_{n+2}+\Delta _{2}\left(
n:n_{1},...,n_{m}\right) -1)\text{.}
\end{eqnarray*}
\end{lemma}

\begin{proof}
We compute the first formula:

\begin{eqnarray*}
\left( c,d\right) \left( \kappa ^{\ast }\right) ^{n_{1}}\sigma ^{\ast
}...\left( \kappa ^{\ast }\right) ^{n_{m}}\sigma ^{\ast } &=&\left(
c,d\right) \left( \sigma ^{\ast }\right) ^{\left( \lambda ^{\ast }\right)
^{-n_{1}}}\left( \sigma ^{\ast }\right) ^{\left( \lambda ^{\ast }\right)
^{-n_{1}-n_{2}}}...\left( \sigma ^{\ast }\right) ^{\left( \lambda ^{\ast
}\right) ^{-n}}\left( \kappa ^{\ast }\right) ^{n}, \\
&=&(c+a_{-n_{1}-1}+a_{-n_{1}-n_{2}-1}+...+a_{-n-1}, \\
&&d+a_{-n_{1}}+a_{-n_{1}-n_{2}}+...+a_{-n})\kappa ^{n} \\
&=&(a_{n-1}\left( c+a_{-n_{1}-1}+a_{-n_{1}-n_{2}-1}+...+a_{-n-1}\right) \\
&&+a_{n}\left( d+a_{-n_{1}}+a_{-n_{1}-n_{2}}+...+a_{-n}\right) +a_{n+1}-1, \\
&&a_{n}\left( c+a_{-n_{1}-1}+a_{-n_{1}-n_{2}-1}+...+a_{-n-1}\right) \\
&&+a_{n+1}\left( d+a_{-n_{1}}+a_{-n_{1}-n_{2}}+...+a_{-n}\right) \\
&&+a_{n+2}-1) \\
&=&... \\
&=&(a_{n-1}c+a_{n}d+a_{n+1}+ \\
&&\left( a_{n-n_{1}-1}+a_{n-\left( n_{1}+n_{2}\right) -1}+...+a_{n-\left(
n_{1}+n_{2}+...+n_{m-1}\right) -1}\right) , \\
&&a_{n}c+a_{n+1}d+a_{n+2}+ \\
&&\left( a_{n-n_{1}}+a_{n-\left( n_{1}+n_{2}\right) }+...+a_{n-\left(
n_{1}+n_{2}+...+n_{m-1}\right) }\right) -1)\text{.}
\end{eqnarray*}%
The other formulae can be derived similarly.
\end{proof}

On substituting $c=0,d=1$, in the above, we get the following set of
possible types $\left( c_{n},d_{n}\right) $ and the corresponding lengths of
elements of $W$:

\begin{eqnarray*}
\text{(i) }\left( 0,1\right) \left( \kappa ^{\ast }\right) ^{n_{1}}\sigma
^{\ast }...\left( \kappa ^{\ast }\right) ^{n_{m}}\sigma ^{\ast }
&=&(a_{n+2}+\Delta _{1}\left( n:n_{1},...,n_{m}\right) , \\
&&a_{n+3}-1+\Delta _{2}\left( n:n_{1},...,n_{m}\right) ) \\
&&\text{ of length }\left( a_{n+4}-1\right) +\Delta \left(
n:n_{1},...,n_{m}\right) \text{;} \\
\text{(ii) }\left( 0,1\right) \left( \kappa ^{\ast }\right) ^{n_{1}}\sigma
^{\ast }...\left( \kappa ^{\ast }\right) ^{n_{m}} &=&(a_{n+2}-1+\Delta
_{1}\left( n:n_{1},...,n_{m}\right) , \\
&&a_{n+3}-1+\Delta _{2}\left( n:n_{1},...,n_{m}\right) ) \\
&&\text{of length }\left( a_{n+4}-2\right) +\Delta \left(
n:n_{1},...,n_{m}\right) \text{;} \\
\text{(iii) }\left( 0,1\right) \sigma ^{\ast }\left( \kappa ^{\ast }\right)
^{n_{1}}\sigma ^{\ast }...\left( \kappa ^{\ast }\right) ^{n_{m}}\sigma
^{\ast } &=&(2a_{n+1}+\Delta _{1}\left( n:n_{1},...,n_{m}\right) , \\
&&2a_{n+2}-1+\Delta _{2}\left( n:n_{1},...,n_{m}\right) ) \\
&&\text{of length }\left( 2a_{n+3}-1\right) +\Delta \left(
n:n_{1},...,n_{m}\right) \text{;} \\
\text{(iv) }\left( 0,1\right) \sigma ^{\ast }\left( \kappa ^{\ast }\right)
^{n_{1}}\sigma ^{\ast }...\left( \kappa ^{\ast }\right) ^{n_{m}}
&=&(2a_{n+1}-1+\Delta _{1}\left( n:n_{1},...,n_{m}\right) , \\
&&2a_{n+2}-1+\Delta _{2}\left( n:n_{1},...,n_{m}\right) ) \\
&&\text{of length }\left( 2a_{n+3}-2\right) +\Delta \left(
n:n_{1},...,n_{m}\right) \text{.}
\end{eqnarray*}

\subsubsection{Lengths of words and Fibonacci lengths}

Suppose an integer $\Delta $ is in the range $\left( 0,a_{n+2}-2\right) $.
Then, there exists at least one partition $n_{1},...,n_{m}$ of $n$ such that 
$\Delta =\Delta \left( n:n_{1},...,n_{m}\right) $: let $a_{k_{1}}\leq \Delta
<a_{k_{1}+1}$; then $k_{1}\leq n+1$. Write $\Delta $ in base Fibonacci: $%
\Delta =a_{k_{1}}+a_{k_{2}}+...+a_{k_{r}}$ with $r\geq 1$ and $%
k_{1}>k_{2}>...>k_{r}\geq 1$. We call $r$ the Fibonacci length of $\Delta $.
Then $r<n$ : for the maximum Fibonacci length of $\Delta $ is $k_{1}$, in
which case, $\Delta =a_{k_{1}}+a_{k_{1}-1}+...+a_{1}=a_{k_{1}+2}-1<a_{n+2}-2$
and therefore, $k_{1}<n$. Now, the partition $\left( n_{1},...,n_{m}\right) $
of $n$ is determined by%
\begin{eqnarray*}
n_{1} &=&n-k_{1}+1, \\
n_{i} &=&n-\left( n_{1}+n_{2}+...+n_{i-1}+k_{i}\right) +1\text{ for }1<i<r,
\\
n_{r} &=&n-\left( n_{1}+n_{2}+...+n_{r-1}+k_{r}\right) +1, \\
n_{r+1} &=&n-\left( n_{1}+...+n_{r}\right) \text{.}
\end{eqnarray*}

Let $w\in W$ have length $l\left( w\right) $. Let $l\left( w\right) $ lie in
the range $[a_{n+4}-2,a_{n+5}-2)$ for $n\geq 0$; then $n$ is unique. The
length $l\left( w\right) $ has four possible forms%
\begin{eqnarray*}
&&\left( a_{n+4}-1\right) +\Delta \left( n:n_{1},...,n_{m}\right) ,\left(
a_{n+4}-2\right) +\Delta \left( n:n_{1},...,n_{m}\right) , \\
&&\left( 2a_{n+3}-1\right) +\Delta \left( n:n_{1},...,n_{m}\right) ,\left(
2a_{n+3}-2\right) +\Delta \left( n:n_{1},...,n_{m}\right) \text{.}
\end{eqnarray*}%
Each unordered partition $n_{1},...,n_{m}$ of $n$ produces a $\Delta \left(
n:n_{1},...,n_{m}\right) $ which determines $\Delta _{1}\left(
n:n_{1},...,n_{m}\right) ,\Delta _{2}\left( n:n_{1},...,n_{m}\right) $ and
thus produces $4$ different words with lengths in the given range. On the
other hand, given an integer $l$ in the given range, then at least one of
the values $l-\left( a_{n+4}-2\right) ,l-\left( a_{n+4}-1\right) ,l-\left(
2a_{n+3}-2\right) ,l-\left( 2a_{n+3}-1\right) $, call it $\Delta $, is in
the range $[0,a_{n+2}-2]$ and therefore there exists a partition $%
n_{1},...,n_{m}$ such that $\Delta =\Delta \left( n:n_{1},...,n_{m}\right) $%
. There are $2^{n-1}$ unordered partitions of $n$. Therefore, the number of
different monomials $w\in W$ with lengths in the range $%
[a_{n+4}-2,a_{n+5}-2) $ is $4.2^{n-1}=2^{n+1}$.

Now we compute, for $n\geq 4$ 
\begin{eqnarray*}
f\left( a_{n+1}-3\right) &=&f\left( a_{n+1}-4\right) +g\left(
a_{n+1}-3\right) \\
&=&f\left( a_{n+1}-3-i\right) +g\left( a_{n+1}-3-i+1\right) +...+g\left(
a_{n+1}-3\right) \\
&=&f\left( a_{n+1}-3-a_{n-1}\right) +g\left( a_{n+1}-3-a_{n-1}+1\right)
+...+g\left( a_{n+1}-3\right) \\
&=&f\left( a_{n}-3\right) +\left( g\left( a_{n}-2\right) +...+g\left(
a_{n+1}-3\right) \right) \\
&=&f\left( a_{n}-3\right) +2^{n-3} \\
&=&f\left( a_{4}-3\right) +2^{2}+...+2^{n-3} \\
&=&2^{n-2}-1\text{.}
\end{eqnarray*}%
Thus, as $g\left( a_{n+1}-2\right) =1$, we obtain the desired 
\begin{equation*}
f\left( a_{n+1}-2\right) =2^{n-2}\text{.}
\end{equation*}
Hence, by Lemma 9, we have proved the growth type of $f$ equal to $n^{\frac{%
\log \left( 2\right) }{\log \left( \alpha \right) }}$, that of $F$ equal to $%
n^{1+\frac{\log \left( 2\right) }{\log \left( \alpha \right) }}$ and thence
the Gelfand- Kirillov dimension of $R$ equal to $1+\frac{\log \left(
2\right) }{\log \left( \alpha \right) }$.

\section{Final Comments}

It was shown by M. Vorobets and Y. Vorobets in \cite{voro07} that, in our
notation, the group generated by 
\begin{equation*}
a=\left( 
\begin{array}{cc}
0 & b \\ 
c & 0%
\end{array}%
\right) ,b=\left( 
\begin{array}{cc}
a & 0 \\ 
0 & a%
\end{array}%
\right) ,c=\left( 
\begin{array}{cc}
0 & c \\ 
b & 0%
\end{array}%
\right) \in M\left( \mathbf{T}_{2},\mathbb{Q}\right)
\end{equation*}%
is free of rank $3$. It would be interesting to know whether the algebra $%
\left\langle a,b,c\right\rangle $ is a free group algebra.

With respect to the second ring of this paper, we recall \ a construction of
a $2$-generated ring having its multiplicative semigroup of monomials nil of
degree $3$, which is obtained simply as follows. Consider the free ring $%
F=\left\langle a,b\right\rangle $ in two generators and in it the set $M$ of
all monomials in $a,b$ which contain cubic subwords. The complement of $M$
contains the Morse-Thue sequence (generated by the substitution $%
a\rightarrow ab,b\rightarrow ba$) \cite{mor-hed44}) and is therefore
infinite. The linear closure of $M$ in $F$ is an ideal and $\frac{F}{M}$ is
an infinite ring where each monomial is nil of degree at most $3$. Does the
ring $\frac{F}{M}$ admit a faithful finite-state representation?

In the context of this paper, it is important to decide whether there exist
associative finitely generated infinite nil algebras which are recursive and
finite-state. V. Petrogradsky \cite{petro06} and I. Schestakov-E.Zelmanov 
\cite{shes-zel08} have produced new exciting constructions of nil Lie
algebras in finite characteristic, which are generated by recursively
defined derivations of a polynomial ring. These constructions may prove to
be relevant to the problem for associative algebras.

Bartholdi and Reznykov studied in \cite{bar-rez08} the semigroup generated
by the recursive matrices 
\begin{equation*}
s=\left( 
\begin{array}{cc}
0 & 1 \\ 
1 & 0%
\end{array}%
\right) ,t=\left( 
\begin{array}{cc}
s & 0 \\ 
t & 0%
\end{array}%
\right)
\end{equation*}%
with $s_{\phi }=t_{\phi }=1$ and showed that this semigroup satisfies the
identity $w^{6}=w^{4}$.

We thank R. Grigorchuk for his insistence on knowing the growth of $R_{2}$
and thank Ricardo Nunes de Oliveira for computer computations which were
helpful in answering the question. We note that the types of growth for our
ring $R_{2}$, of the Petrogradsky Lie algebra and that of Bartholdi-Reznykov
semigroup are all equal.

\end{document}